\documentclass[12pt]{amsart}
\usepackage{amsmath, amsfonts, amsthm}
\usepackage{amssymb}
\usepackage{mathrsfs}
\usepackage[a4paper, margin=1in]{geometry}
\usepackage{tikz-cd}
\usepackage{enumitem}
\usepackage{quiver}
\usepackage{array,makecell}
\usepackage[table]{xcolor}
\usepackage{diagbox}

\usepackage{hyperref}
\hypersetup{
    colorlinks=true,
    linkcolor=blue,
    filecolor=magenta,      
    urlcolor=cyan,
    citecolor=red
}
\newcommand{\Mod}[1]{\ (\mathrm{mod}\ #1)}

\newtheorem{theorem}{Theorem}[section]

\newtheorem{lemma}[theorem]{Lemma}

\theoremstyle{remark}
\newtheorem{remark}[theorem]{Remark}

\theoremstyle{definition}
\newtheorem{definition}[theorem]{Definition}

\newcommand{\qbf}{\mathbf{q}}
\newcommand{\rbf}{\mathbf{r}}
\newcommand{\dbf}{\mathbf{d}}

\newcommand{\Mcal}{\mathcal{M}}

\newcommand{\Ocal}{\mathcal{O}}

\newcommand{\Xcal}{\mathcal{X}}

\newcommand{\Scal}{\mathcal{S}}

\newcommand{\Rcal}{\mathcal{R}}
\newcommand{\Zcal}{\mathcal{Z}}

\newcommand{\Tcal}{\mathcal{T}}

\newcommand{\Ycal}{\mathcal{Y}}

\newcommand{\PP}{\mathbb{P}}
\renewcommand{\AA}{\mathbb{A}}

\newcommand{\ZZ}{\mathbb{Z}}

\newcommand{\GG}{\mathbb{G}}

\newcommand{\Sfrak}{\mathfrak{S}}
\newcommand{\gfrak}{\mathfrak{g}}

\newcommand{\gitquot}{/ \! \! /}

\newcommand{\length}{\mathrm{length}}

\renewcommand{\H}{\mathrm{H}}

\newcommand{\Pic}{\mathrm{Pic}}

\newcommand{\Aut}{\mathrm{Aut}}

\newcommand{\PGL}{\mathrm{PGL}}

\newcommand{\pr}{\mathrm{pr}}

\begin{document}
\author{Bogdan Carasca \& Vlad Robu}
\date{\today}
\address{Humboldt Universität zu Berlin, Unter den Linden 6, 10117 Berlin}
\email{bogdan.carasca@hu-berlin.de}
\email{vlad.nicolae.robu@hu-berlin.de}

\begin{abstract}
    We prove that all components of the meromorphic differential strata are unirational in genus $2\le g\le 5$ and with relatively low number of poles. In particular, it follows that the moduli space of pointed even spin curves $\Scal_{5;1}^+=\Scal_5^+\times_{\Mcal_5}\Mcal_{5,1}$ is unirational. 
\end{abstract}

\title[Unirationality of Meromorphic Differential Strata]{Unirationality of Meromorphic Differential Strata}
\maketitle

\let\thefootnote\relax

\section{Introduction}

Let $g\ge 2$ be an integer and let $\mu=(\mu_1,\dots,\mu_n)$ be a partition with nonzero entries of $2g-2$. The moduli space parametrising smooth algebraic curves of genus $g$ equipped with a canonical divisor having marked zeroes and/or poles of orders prescribed by $\mu$ is denoted by
$$\PP\Omega_g(\mu):=\left\{\left[C,p_1, \dots, p_n\right]\in \Mcal_{g,n}\colon\Ocal_C\left(\sum\limits_{i=1}^n\mu_ip_i\right)\cong \omega_C\right\}.$$
Depending on whether $\mu$ contains only nonnegative entries or not, we refer to $\PP\Omega_g(\mu)$ as a \textit{holomorphic} or a \textit{meromorphic} stratum of differentials. In the holomorphic case, it is known that $\dim \PP\Omega_g(\mu)=2g-2+n$ (see e.g. \cite{polishchuk2006}), whereas in the meromorphic case we have that $\dim \PP\Omega_g(\mu)=2g-3+n$, thanks to \cite{boissy2015connCmpDiffStrata}.
\vspace{0.5em}

The birational geometry of the moduli spaces of stable curves $\overline{\Mcal}_g$ is a central topic in the subject. It was known since Severi \cite{severi} that $\overline{\Mcal}_g$ is unirational for $g\le 10$. Much later, the celebrated papers \cite{harrisMumford} and \cite{harris1987kodairaDimMg} showed that $\overline{\Mcal}_g$ is of general type, when $g\ge 24$. More recently, the moduli spaces  $\overline{\Mcal}_{22}$ and $\overline{\Mcal}_{23}$ were also settled to have maximal Kodaira dimension in \cite{farkasjensenpayne25}. Furthermore, it is now known that $\overline{\Mcal}_g$ is unirational for $g\le 14$, due to \cite{verra14}, that $\overline{\Mcal}_{15}$ is rationally connected by \cite{brunoverra15}, and the recent paper \cite{farkasverra16} shows that $\overline{\Mcal}_{16}$ is uniruled. 

\vspace{0.5em}
In light of these developments, the study of the birational geometry of strata of differentials is only natural to pursue. The complete Kodaira classification of the the closely related odd spin moduli spaces $\Scal_g^-$ has been achieved in \cite{FarkasVerraOddSpin}. These spaces are birational to the quotient of $\PP\Omega_g(2^{g-1})^-$ by the action of the symmetric group $\Sfrak_{g-1}$. More generally, in the case of holomorphic strata with arbitrary positive multiplicities, \cite{barros2018} and \cite{bud2021} offer uniform results that cover certain low genus cases. Concretely, they show that all components of the holomorphic stratum $\PP\Omega_g(\mu)$ are uniruled in the following ranges:
\begin{itemize}
    \item when $2\le g \le 9$ and all possible holomorphic partitions, 
    \item when $g=10$ and $11\le n<18$, and  
    \item when $g=11$, $n\ge 10$ and $\mu\ne \left(2^{10}\right)$.
\end{itemize}

At the same time, the odd spin component of $\PP\Omega_g(2g-2)$ is of general type whenever $g\ge 13$, as of \cite{chencostantinimoeller2024}, and the even spin component of the same space is of general type whenever $g=31$ or $g\ge 33$, thanks to \cite{budchenmoeller2024}. 

\vspace{1em}
Surprisingly, nothing seems to be known about the birational geometry of meromorphic strata. In this paper, we present a comprehensive study of these questions in the meromorphic case for low genus.

\vspace{1em}
To make things transparent, we consistently separate the positive and negative parts of the partition. Precisely, $\mu=(\mu_1,\dots,\mu_n)$ and $\nu=(\nu_1,\dots,\nu_r)$ stand for decreasing tuples of positive integers satisfying
\begin{align}\label{eq: sum is 2g-2}
    \sum\limits_{i=1}^n\mu_i-\sum\limits_{j=1}^r\nu_j=2g-2,
\end{align}
with $\mu$ and $\nu$ indicating the zeroes, respectively the poles of the abelian differential. We write 
$$\PP\Omega_g(\mu,\nu)=\left\{\left[C,(p_i)_{i=1}^n,(q_j)_{j=1}^r\right]\in \Mcal_{g,n+r}:\Ocal\left(\sum_{i=1}^n\mu_ip_i-\sum_{j=1}^r\nu_jq_j\right)\cong \omega_C\right\}.$$
We also let $|\mu|$ and $|\nu|$ denote the sum of the entries of $\mu$ and $\nu$, respectively, so in particular $|\mu|-|\nu|=2g-2$ and $|\nu|\ge 2$ for a nonempty stratum $\PP\Omega_g(\mu,\nu)$.

\vspace{1em}
For $g\ge 2$, it was shown in \cite{boissy2015connCmpDiffStrata} that $\PP\Omega_g(\mu,\nu)$ has at most $3$ irreducible components, and their exact number can be read directly from the signatures $\mu$ and $\nu$. To be precise, we call a partition \textit{even} if all its entries are even and \textit{odd} otherwise. Thus we say:

\begin{itemize}
    \item a partition $(\mu,\nu)$ is of \textit{even type}, if either both $\mu$ and $\nu$ are even, or if $\mu$ is even and $\nu=(1,1)$, and 
    \item a partition $(\mu,\nu)$ is of \textit{hyperelliptic type}, if both $\mu$, $\nu$ are of the form $(a,a)$ or $(2a)$.
\end{itemize}

\vskip 0.5em
\noindent The number of components is equal to 
\begin{enumerate}[label=\textbf{(\alph*)}]
    \item\label{not even, not hyper} one, if $(\mu,\nu)$ is neither of even type, nor of hyperelliptic type, 
    \item\label{even, not hyper} two, if $(\mu,\nu)$ is of even type, but not of hyperelliptic type, 
    \item\label{not even, hyper} two, if $(\mu,\nu)$ is of hyperelliptic type, but not of even type, and 
    \item\label{even, hyper} three, if $(\mu,\nu)$ is of both even type and hyperelliptic type. 
\end{enumerate}

\vskip 1em
\noindent There is one single exception, namely in genus $g=2$ with $|\nu|=2$, when there are only two components, for partitions $(\mu,\nu)$ that are of both even and hyperelliptic type. In the cases \ref{not even, hyper} and \ref{even, hyper}, as well as in this exception, there is one anti--invariant hyperelliptic component, which is already known to be unirational (see e.g. \cite[Section 2.3]{chencostantinimoeller2024}). Consequently, in order to to avoid repetition, in the text we simply exclude this special component. If $(\mu,\nu)$ is of even type, then the two components $\PP\Omega_g(\mu,\nu)^+$ and $\PP\Omega_g(\mu,\nu)^-$ are indexed by the parity of the spin structure, that is 
\begin{itemize}
    \item the parity of $h^0\left(\sum\limits_{i=1}^n\mu_i/2\cdot p_i-\sum\limits_{j=1}^r\nu_j/2\cdot q_j\right)$ when $\mu$ and $\nu$ are even, or
    \item the parity of $h^0\left(\sum\limits_{i=1}^n\mu_i/2\cdot p_i\right)$ when $\mu$ is even and $\nu=(1,1)$. 
\end{itemize} 

\vskip 1em
Our main result is summarized in the next theorem. The word "always" appearing in the table below has the following precise meaning: given the genus $g$ and the weight $k$ of the negative part of the partition, unirationality of the stratum $\PP\Omega_g(\mu, \nu)$ holds for every $\nu$ with $|\nu| = k$ and $\mu$ subject to condition (\ref{eq: sum is 2g-2}).

\vskip 2em
\begin{theorem}\label{thm: unirat diff strata in genus 2}\label{thm: unirat diff strata in genus 3}\label{thm: unirat diff strata in genus 4}\label{thm: unirat diff strata in genus 5}
All components of the meromorphic differential stratum $\PP\Omega_g(\mu,\nu)$ are unirational whenever $g$, $\mu$ and $\nu$ are as in the table below:
\vspace{0.5em}
\begin{center}
\renewcommand{\arraystretch}{1.5}
\setlength{\tabcolsep}{5pt}
\begin{tabular}{|c|c|c|c|c|}
\hline
\diagbox[dir=NW,width=2.5cm,height=0.8cm]{}{}
  & $g=2$ & $g=3$ & $g=4$ & $g=5$ \\
\hline
$|\nu|=2$
  & always & always & always & $n\le 9$ \\
\hline
$|\nu|=3$
  & always & always & $n\le 8$ & \\
\hline
$|\nu|=4$
  & always & always & & \\
\hline
$|\nu|=5$
  & always & & & \\
\hline
$|\nu|=6$
  & \makecell{$\nu_{r-2}\ge 2$ and $n\le 7$ if $\nu$ odd}
  & & & \\
\hline
$|\nu|=7$
  & \makecell{$\nu=(2,2,2,1),\ (4,2,1),$\\
              $(3,2,2),\ (5,2)$ and $n\le 8$}
  & & & \\
\hline
\end{tabular}
\end{center}
\end{theorem}
\vspace{0.4em}

Theorem \ref{thm: unirat diff strata in genus 5} has interesting applications to the moduli spaces
$$
\Scal_g^+=\left\{\left[C,\eta\right]:\left[C\right]\in\Mcal_g\textnormal{ and }\eta^{\otimes 2}\cong \omega_C\textnormal{, with }h^0(C,\eta)\equiv 0\Mod{2}\right\}
$$ 
of even spin curves of genus $g$. In \cite{dolgachevkanev1993} and \cite{takagiRatS4plus}, it is shown that $\Scal_3^+$ and $\Scal_4^+$ are rational varieties, respectively. Compared to the existing methods leading to the rationality of $\Scal_3^+$ and $\Scal_4^+$, the merit of our proofs lies in the simplicity of the models. To complete the picture of the Kodaira classification of $\Scal_g^+$, we recall that Farkas and Verra establish that $\Scal_g^+$ is uniruled for $g\le 7$, of Calabi--Yau type for $g=8$ and of general type for $g\ge 9$ (\cite{farkasEvenSpin} and \cite{farkasverra2012}), so the problem of the Kodaira dimension of $\Scal_g^+$ is settled for all genera. In this context, it is interesting to see that Theorem \ref{thm: unirat diff strata in genus 2} leads to the following improvement in genus $5$:

\begin{theorem}\label{thm: unirat S5 even}
The moduli space of pointed even spin curves $\Scal_{5;1}^+:=\Scal_5^+\times_{\Mcal_5}\Mcal_{5,1}$ is unirational. In particular, $\Scal_5^+$ is unirational.
\end{theorem}

Indeed, for the signatures $\mu=(2,2,2,2,2)$ and $\nu=(2)$, there is a dominant map 
$$\PP\Omega_5(\mu,\nu)^+\dashrightarrow \Scal_{5;1}^+\textnormal{ defined by }\left[C,p_1,\dots,p_5,q\right]\mapsto \left[C,\Ocal_C(p_1+\dots+p_5-q),q\right],$$
so the unirationality of $\PP\Omega_5(\mu,\nu)^+$ immediately implies the unirationality of $\Scal_{5;1}^+$.
\vspace{0.5em}

Returning to Theorem \ref{thm: unirat diff strata in genus 2}, all its cases differ by virtue of the distinct geometry of each generic point of $\PP\Omega_g(\mu,\nu)$. There is, nonetheless, a common thread to all our proofs, which we now summarise in the sequence of steps below.

\vspace{0.5em}
\textbf{(1)} We embed the curve $C$ underlying a general point in a component of $\PP\Omega_g(\mu,\nu)$ in a rational surface $S$. 

\vspace{0.5em}
\textbf{(2)} We show that $S$ contains curves $D$ such that $C\cdot D=\sum\limits_{j=1}^r\nu_jq_j+\dbf$ for some effective divisor $\dbf$.

\vspace{0.5em}
\textbf{(3)} Write $\omega_C=\Ocal_C(C+K_S)$ using adjunction, and then we find the divisor 
$$\sum\limits_{i=1}^n\mu_ip_i+\dbf\in |\Ocal_C(C+K_S+D)|.$$
Here, the restriction map $\H^0(S,\Ocal_S(C+K_S+D))\to \H^0(C,\Ocal_C(C+K_S+D))$ is surjective, so there is a curve $F \in |C + K_S + D|$ satisfying $C\cdot F=\sum\limits_{i=1}^n\mu_ip_i+\dbf$. 

\vspace{0.5em}
\textbf{(4)} The data $\left(C,D,F,(p_i)_{i=1}^n,(q_j)_{j=1}^r,\dbf\right)$ together with the above intersection relations are parametrized by an incidence variety $\Xcal_g(\mu,\nu)$ and we obtain a dominant rational map $\Xcal_g(\mu, \nu) \dashrightarrow \PP\Omega_g(\mu, \nu)$.
\vspace{0.5em}

\textbf{(5)} The task is then to show the rationality of $\Xcal_g(\mu, \nu)$. To that end, we consider the variety $\Ycal_g(\mu, \nu)$ parametrising the data $\left(D,F,(p_i)_{i=1}^n,(q_j)_{j=1}^r,\dbf\right)$ appearing in the previous step. The projection
\[
\pi : \Xcal_g(\mu, \nu) \dashrightarrow \Ycal_g(\mu, \nu)
\]
is then shown to be a projective bundle.

\vskip 0.5em
\textbf{(6)} In this last step, the rationality of $\Ycal_g(\mu, \nu)$ is proven.

\vspace{0.4em}






\begin{remark}
In the context of the steps \textbf{(1)}--\textbf{(6)} described above, the number of (nonhyperelliptic) components always reflects in the structure of the varieties $\Xcal_g(\mu,\nu)$ and $\Ycal_g(\mu,\nu)$. If, say, there are $2$ components, then there exist two forgetful maps 
$$\pi^+:\Xcal_g(\mu,\nu)^+\to \Ycal_g(\mu,\nu)^+\textnormal{ and }\pi^-:\Xcal_g(\mu,\nu)^-\to \Ycal_g(\mu,\nu)^-$$
corresponding to $\PP\Omega_g(\mu,\nu)^+$ and $\PP\Omega_g(\mu,\nu)^-$ respectively, and both these maps are shown to be projective bundles over dense opens in their respective target. 
\end{remark}

\subsection*{Acknowledgements} We thank our advisor Gavril Farkas for his constant guidance and support. We also thank Marian Aprodu, Riccardo Redigolo and Andrei Stoenică for useful conversations on the subject. 

\vskip 0.5em
This work was funded by the Deutsche Forschungsgemeinschaft (DFG,
German Research Foundation) - RTG2965 - Project number 512730679. This work was also funded by the DFG under Germany's Excellence Strategy - The Berlin Mathematics Research Center MATH+.

\subsection*{AI statement} This paper was written in the style of the pre--AI era. As such, we made no use of the existing AI tools in its development.

\section{The Independence Result}






As mentioned in the Introduction, our strategy is to embed the curve $C$ underlying the general point of $\PP\Omega_g(\mu, \nu)$ in a rational surface and thus obtain incidence varieties $\Xcal_g(\mu, \nu) / \Ycal_g(\mu, \nu)$. One of the key points of the proof is then to show that $\Xcal_g(\mu, \nu)$ is a projective bundle over a dense open in $\Ycal_g(\mu, \nu)$. To that end, the following result will prove particularly useful.

\begin{theorem}\label{big indep lemma}
Let $n\in \ZZ_{\ge 1}$ be an integer and let $S$ be a smooth rational surface. Fix an ample line bundle $L$ and a line bundle $M$ on $S$ such that $\omega_S=-nL$, and suppose there is a smooth irreducible curve $C$ in the linear system $|(n+1)L-M|$, where
\[
M\cdot L+2g(C)-2>0.
\]
If $Z$ is a zero dimensional subscheme of $C$, then $Z$ imposes precisely
\[
\length(Z)-h^0(C,Z-nK_C-nM\vert_C)
\]
independent conditions on the linear system $|(n+1)L-M|$.
\end{theorem}

The statement is written in full generality, but here we use it for the values $n\le 3$. Specifically, 

\textbf{(a)} $n=3$ appears when $S=\PP^2$,

\textbf{(b)} $n=2$ appears when $S=\PP^1\times \PP^1$, and

\textbf{(c)} $n=1$ appears when $S$ is a del Pezzo surface.

\begin{proof}
The number of conditions imposed by $Z$ on the linear system $|(n+1)L-M|$ is the rank of the restriction map $\textnormal{res}:\H^0(S,(n+1)L-M)\to \H^0(Z,(n+1)L-M)$. We notice that this map factors as
$$\H^0(S,(n+1)L-M)\stackrel{\textnormal{res}_1}{\longrightarrow} \H^0(C,(n+1)L-M)\stackrel{\textnormal{res}_2}{\longrightarrow} \H^0(Z,(n+1)L-M).$$
Now, the first restriction map is surjective. Indeed, we have the short exact sequence of sheaves 
$$0\to \Ocal_S\to \Ocal_S((n+1)L-M)\to \Ocal_C((n+1)L-M)\to 0$$
and $\H^0(S,(n+1)L-M)\stackrel{\textnormal{res}_1}{\longrightarrow} \H^0(C,(n+1)L-M)\to \H^1(S,\Ocal_S)=0$, where the last equality comes from the fact that $S$ is regular. Consequently, the maps $\textnormal{res}$ and $\textnormal{res}_2$ have the same rank. Furthermore, $\omega_C=\Ocal_C(L-M)$ by adjunction and, if we set $A=\Ocal_C(M)$, we are simply interested in the rank of the map 
$$\textnormal{res}_2:\H^0\left(C,\omega_C^{n+1}\otimes A^{n}\right)\to \H^0\left(Z,\omega_C^{n+1}\otimes A^{n}\right).$$
Finally, passing to cohomology the short exact sequence 
$$0\to \Ocal_C(-Z)\left(\omega_C^{n+1}\otimes A^{n}\right)\to \Ocal_C\left(\omega_C^{n+1}\otimes A^{n}\right)\to \Ocal_Z\left(\omega_C^{n+1}\otimes A^{n}\right)\to 0$$
gives 
$$\H^0\left(C,\Ocal_C\left(\omega_C^{n+1}\otimes A^{n}\right)\right)\stackrel{\textnormal{res}_2}{\longrightarrow} \H^0\left(Z,\Ocal_Z\left(\omega_C^{n+1}\otimes A^{n}\right)\right)\twoheadrightarrow \H^1\left(\Ocal_C(-Z)\left(\omega_C^{n+1}\otimes A^{n}\right)\right),$$
where the last map is a surjection, since we have the vanishing of the next term in the sequence, $\H^1(C,\Ocal_C(\omega_C^{n+1}\otimes A^n))=0$. Here, we used $\deg(\omega_C\otimes A)>0$. It follows that the rank in question is $$\length(Z)-h^1\left(\Ocal_C(-Z)\left(\omega_C^{n+1}\otimes A^{n}\right)\right)=\length(Z)-h^0\left(C,\omega_C^{-n}\otimes A^{-n}(Z)\right).$$ 
The Theorem is now proven. 
\end{proof}

\section{Meromorphic Differential Strata in Genus \texorpdfstring{$5$}{}}

The case of genus $5$ is quite instructive for the methods of the entire paper. This section is concerned with the proof of Theorem \ref{thm: unirat diff strata in genus 5} in this case, which states that all components of the space
$$\PP\Omega_5(\mu,\nu)=\left\{[C,(p_i)_{i=1}^n,q^{\pm}]\in\Mcal_{5,n+r}:\Ocal_C\left(\sum_{i=1}^n\mu_ip_i-q^{+}-q^{-}\right)\cong \omega_C\right\}$$
are unirational, whenever $n\le 9$. Here, the situation $\nu=(2)$ corresponds to $q^+=q^-$, while $\nu=(1,1)$ corresponds to $q^+\ne q^-$. When $C$ is a nontrigonal curve of genus $5$, the canonical bundle presents it as the base locus of a net $\Lambda$ of quadrics in $\PP^4$.

\vspace{0.5em}
In order to obtain a surface $S$ on which $C$ lies, we choose a pencil in $\Lambda$, hence embedding $C$ in a smooth del Pezzo surface $S$ of degree $4$ such that $C\in |-2K_S|$. We realise $S$ as the blow-up of $\PP^2$ at $5$ points $q_1,\dots,q_5$ and call $E_{q_i}$ the exceptional divisor above $q_i$, for $1\le i\le 5$. Since $(-2K_S)\cdot E_{q_5}=2$, we have our candidate: given the two poles $q^{\pm}$, we would like to have $C\cdot E_{q_1}=q^++q^{-}$.

\vspace{0.5em}
Now blow everything down to $\PP^2$. This means that we are interested in plane sextics $\Gamma$ that are singular at $4$ points $q_2,\dots,q_5$ and have a node/cusp at $q_1$, whose tangent cone indicates the points $q^+$ and $q^-$ on $C$. The two lemmas below create precisely this environment. 
\vspace{0.5em}

\begin{lemma}\label{points exchanged by bielliptic struct}
Let $[C]\in \Mcal_{5}$ be a nontrigonal curve, $q^+,q^-$ two points on $C$ (not necessarily distinct) and consider the induced map $\sigma:C\to \Gamma \subset \PP^2$ given by $L=\omega_C(-q^+-q^-)$. 

\vskip 0.5em
Then, $\sigma$ is a $2:1$ cover (in which case $\Gamma$ is an elliptic plane cubic) if and only if $q^+$ and $q^-$ are exchanged by the corresponding bielliptic structure. 
\end{lemma}
\begin{proof}
Any bielliptic structure on the curve $C$ corresponds to presenting the canonical model of $C$ as the complete intersection of a cone over an elliptic normal curve $E\subset \PP^3$ and a quadric in $\PP^4$ (see e.g. \cite{bardellidelcentina1990}). If $\sigma$ is a $2:1$ cover, then the map $\sigma:C\to \Gamma\subset \PP^2$ factors as 
\[\begin{tikzcd}
	C & E & \Gamma
	\arrow["\iota", from=1-1, to=1-2]
	\arrow["\sigma"', curve={height=12pt}, from=1-1, to=1-3]
	\arrow["\pr", from=1-2, to=1-3]
\end{tikzcd}\]
where $\iota$ is the corresponding bielliptic structure, described as above, and $\pr:E\to \Gamma$ is the projection through a point $x\in E$ on $\PP^2$. Then $q^+$ and $q^-$ are precisely the points in the fibre of $\iota$ over $x$, i.e. $\iota$ involutes $q^+$ and $q^-$. Conversely, if $q^+$ and $q^-$ are exchanged by a bielliptic structure $\iota:C\to E\subset \PP^3$, projecting through their image $x\in E$ to $\PP^2$ gives back the map $\sigma:C\to \Gamma$ induced by $L=\omega_C(-q^+-q^-)$. 
\end{proof}

Before moving to the second lemma, we recall that, for a nonhyperelliptic smooth curve $[C]\in \Mcal_5$, there is an isomorphism $C_2\cong G^2_6(C)$, given by $q^++q^-\mapsto \omega_C(-q^+-q^-)$. 

\begin{lemma}\label{5 systems having a fixed cusp}
Let $[C]\in \Mcal_5$ be a nontrigonal curve and let $q^+,q^-$ be two points on $C$, not necessarily distinct, such that they are not exchanged by a bielliptic structure on $C$, if there are any. 
\vspace{0.5em}

Then, there are five $\gfrak^2_6$'s on $C$ such that the image curve of their induced map $\sigma:C\to \Gamma\subset \PP^2$ is a birational image of $C$, having a node below $q^+$ and $q^-$. If $q:=q^+=q^-$, then we mean the image curve $\Gamma$ has a cusp at $q$. 
\end{lemma}
\begin{proof}
Let $x,y\in C$ be two points defining a $\gfrak^2_6$ via $x+y\mapsto L:=\omega_C(-x-y)$ such that $\sigma_L:C\stackrel{|L|}{\longrightarrow}\PP^2$ is birational onto its image. Then $q^+,q^-$ are in the same fibre of $\nu$ if and only if $h^0(C,\omega_C(-x-y-q^+-q^-))\ge 2$, which is the same as requiring $x,y$ be in the same fibre of the map induced by $M:=\omega_C(-q^+-q^-)$. In other words, since Lemma \ref{points exchanged by bielliptic struct} tells us that $\sigma_M:C\stackrel{|M|}{\longrightarrow}\PP^2$ is birational onto its image, we only need to count the number of singular points of $\sigma_M(C)$. Indeed, notice that 
\vspace{0.5em}

$\bullet$ $C$ is not trigonal, so that such singular points are just double points, and 

$\bullet$ no double point $(x,y)$ is exchanged by any bielliptic structure $\iota$, as otherwise the condition $h^0(C,\omega_C(-x-y-q^+-q^-))\ge 2$ would force $(q^+,q^-)$ to be exchanged by $\iota$ too.

\vskip 0.5em
Since the image curve $\sigma_M(C)$ has degree $6$ and $g(C)=5$, there are $5$ singular points.  
\end{proof}


\vspace{0.3em}

All preparatory tools are now in place.

\begin{definition}\label{unirat model for merom genus 5}
Consider the incidence variety
$$\Xcal_5(\mu,\nu)=\left\{(\Gamma,E,o,(p_i)_{i=1}^n,q,\qbf,\ell_{\pm}):\begin{matrix}
    E\in |\Ocal_{\PP^2}(3)|, o\in E\textnormal{ and }\qbf\in |\Ocal_E(1)(o)|,\\
    \Gamma\in |\Ocal_{\PP^2}(6)|\textnormal{ nodal along }\qbf\textnormal{ and }q,\\
    \Gamma\cdot E=2\qbf+\mu_1p_1 + \dots + \mu_np_n ,\\
    \textnormal{and }\ell_{\pm}\textnormal{ are the tangent directions at }q
\end{matrix}\right\}.$$
When $\nu=(2)$, the notation $\ell_{\pm}$ stands for a single line and $\Gamma$ is cuspidal at $q$. 
\end{definition}

We are ready to prove the unirationality (of all the components) of $\PP\Omega_5(\mu,\nu)$, when $\nu=(2)$ and $\nu=(1,1)$. Indeed, this is a consequence of the following two lemmas. 

\begin{lemma}\label{dominance lemma in gen 5}
The map $\Xcal_5(\mu,\nu)\dashrightarrow \PP\Omega_5(\mu,\nu)$ given by 
$$\left[\Gamma,E,o,(p_i)_{i=1}^n,q,\qbf,\ell_{\pm}\right]\mapsto \left[C,(p_i)_{i=1}^n,q^+,q^-\right]$$
is well--defined and dominant over each component of the target, where $\sigma:C\to \Gamma$ is the normalisation map and $q^{\pm}$ are the points on $C$ corresponding to the tangent directions $\ell_{\pm}$ to $q$, respectively. 
\vspace{0.5em}

Furthermore, this induces a generic $5:1$ map $\Xcal_5(\mu,\nu)\gitquot\PGL(3)\dashrightarrow \PP\Omega_5(\mu,\nu)$ over each component of the target. 
\end{lemma}
\begin{proof}
Pick a point $\left[\Gamma,E,o,(p_i)_{i=1}^n,q,\qbf,\ell_{\pm}\right]$ in $\Xcal_5(\mu,\nu)$ with $p_g(\Gamma)=5$ and write $\qbf=q_1+\dots+q_4$, for some points $q_j\in \Gamma$. Let $\sigma^{-1}(q_j)=\{q_j^{\pm}\}$, define $\qbf^+=q_1^+\dots+q_4^+$ and $\qbf^-=q_1^-+\dots+q_4^-$, so we can write
$$\omega_C(\qbf^++\qbf^-+q^++q^-)=\sigma^*\omega_\Gamma=\sigma^*\Ocal_\Gamma(3)=\Ocal_C\left(\sum_{i=1}^n\mu_ip_i+\qbf^++\qbf^-\right),$$
whence $\Ocal_C\left(\sum\limits_{i=1}^n\mu_ip_i-q^+-q^-\right)\cong \omega_C$. Now we show the dominance claim. Choose a point $\left[C,(p_i)_{i=1}^n,q^{\pm}\right]$ such that $C$ is not trigonal, nor bielliptic (this generic property is met, as described in Definition \ref{unirat model for merom genus 5}: $5$--nodal/cuspidal plane sextics $\Gamma$ are not trigonal, nor bielliptic). From Lemma \ref{5 systems having a fixed cusp}, it follows that there exist precisely $5$ systems $L\in G^2_6(C)$ such that the induced map $\sigma:C\stackrel{|L|}{\longrightarrow}\Gamma\subset \PP^2$ presents the curve as a plane sextic $\Gamma$ with a node $q$ having $\sigma^{-1}(q)=\{q^+,q^-\}$. If $\nu=(2)$, this means that $\Gamma$ has a cusp at $q$. 
\vspace{0.5em}

The sextic $\Gamma$ has $4$ others singular points $q_1,\dots,q_4$. The relation $\sum\limits_{i=1}^n\mu_ip_i-q^+-q^-\sim  K_C$ implies the relation 
$$\sum\limits_{i=1}^n\mu_ip_i+\sum_{j=1}^4(q_j^++q_j^-)\sim K_C\left(\sum_{j=1}^4(q_j^++q_j^-)+(q^++q^-)\right)=\sigma^*\omega_\Gamma=\sigma^*\Ocal_\Gamma(3).$$
Since the restriction map $\H^0(\PP^2,\Ocal_{\PP^2}(3))\to \H^0(\Gamma,\Ocal_{\Gamma}(3))$ is an isomorphism, we deduce that there exists a unique elliptic curve $E\in |\Ocal_{\PP^2}(3)|$ such that
\[
\Gamma\cdot E=\sum_{i=1}^n\mu_ip_i+2(q_1+\dots+q_4).
\]
Finally, set $\qbf:=q_1+\dots+q_4$ and notice that there exists a unique point $o\in E$ such that $\qbf\in |\Ocal_E(1)(o)|$. The lemma is now proven. 
\end{proof}

We also have a lemma concerned with the rationality of $\Xcal_5(\mu,\nu)$. 

\begin{lemma}\label{rat check Ycal5}
All components of $\Xcal_5(\mu,\nu)$ are rational varieties, whenever $\mu=(\mu_1,\dots,\mu_n)$ is a partition of $10$ with $n\le 9$ and $|\nu|=2$. 
\end{lemma}
\begin{proof}
We introduce the variety 
$$\Ycal_5(\mu)=\left\{(E,o,(p_i)_{i=1}^n,\qbf):\begin{matrix}
    E\in |\Ocal_{\PP^2}(3)|,\qbf\in |\Ocal_E(1)(o)|\textnormal{ and }\\
    p_i\in E\textnormal{ with }\sum\limits_{i=1}^n\mu_ip_i+2o\in |\Ocal_E(4)|
\end{matrix}\right\},$$
as well as the map 
$$\pi:\Xcal_5(\mu,\nu)\dashrightarrow \Ycal_5(\mu)\times \left\{(q,\ell_{\pm}):q\in \PP^2,\textnormal{ }\ell_{\pm}\textnormal{ lines through }q\right\}$$ 
that forgets the curve $\Gamma$. This is well-defined, since the condition $\Gamma\cdot E=\sum\limits_{i=1}^n\mu_ip_i+2\qbf$ forces $\sum\limits_{i=1}^n\mu_ip_i+2\qbf\in |\Ocal_E(6)|$, which in turn gives $\sum\limits_{i=1}^n\mu_ip_i+2o\in |\Ocal_E(4)|$. We check the following three claims, which suffice to conclude that all components of $\Xcal_5(\mu,\nu)$ are rational:
\begin{enumerate}[label=\textbf{(\roman*)}]
    \item\label{indep lemma gen 5} the map $\pi$ is a $\PP^1$-bundle over its image in any of the components of the target, 
    \item\label{rat check for Ycal5} all components of $\Ycal_5(\mu)$ are rational, and
    \item\label{dim count gen 5} $\dim \Xcal_5(\mu,\nu)=15+n+r$ and $\dim \Ycal_5(\mu)=12+n$, so $\pi$ is, in fact, a $\PP^1$-bundle over a dense open in each component of the target. 
\end{enumerate}

The statement \ref{indep lemma gen 5} follows from Theorem \ref{big indep lemma}. To see that, consider the blow-up surface $S:=\textnormal{Bl}_{\qbf,q}(\PP^2)$, and denote by $C\in |6H-2E_\qbf-2E_{q}|$ the strict transform of $\Gamma$. Let the strict transform of $E$ be still denoted by $E\in |3H-E_\qbf|$ and let $q^{\pm}$ be the points on $E_{q}$ corresponding to the directions $\ell_{\pm}$ through $q$. Then the curve $C$ is subject to $C\cdot E=\sum\limits_{i=1}^n\mu_ip_i$ and $C\cdot E_{q}=q^{+}+q^{-}$. In other words, the divisor of conditions is $Z=\sum\limits_{i=1}^n\mu_ip_i+q^++q^-$, and the number of independent conditions imposed on on the system $|6H-2E_\qbf-2E_{q}|\cong \PP^{12}$ is $\length(Z)-h^0(C,Z-K_C)=12-h^0(C,2q^++2q^-)=11$, which settles the first claim. 
\vspace{0.5em}

The dimensionality statement on $\Xcal_5(\mu,\nu)$ in \ref{dim count gen 5} is an immediate consequence of Lemma \ref{dominance lemma in gen 5}. The dimensionality check on $\Ycal_5(\mu)$, as well as \ref{rat check for Ycal5} need to be analysed step by step.
\vspace{0.5em}

\textbf{(a)} \textit{Case of odd $\mu$ and $n\le 9$.} Here, $\PP\Omega_5(\mu,\nu)$ consists of a single component. We expect the same from $\Ycal_5(\mu)$; this follows as a by-product of the following argument. Define the variety 
$$\Zcal_5(\mu)=\left\{\left(E,(p_i)_{i=1}^{n},o\right):E\in |\Ocal_{\PP^2}(3)|\textnormal{, }\sum\limits_{i=1}^n\mu_ip_i+2o\in |\Ocal_E(4)|\right\},$$
so that $\Ycal_5(\mu)\dashrightarrow \Zcal_5(\mu)$ is a $\PP^3$--bundle. It is, thus, sufficient to show that $\Zcal_5(\mu)$ is rational. If $\mu_n=1$, then there is a birational isomorphism 
$$\Zcal_5(\mu)\dashrightarrow \left\{\left(E,(p_i)_{i=1}^{n-1},o\right):E\in |\Ocal_{\PP^2}(3)|\textnormal{, }p_i,o\in E\right\},$$
since $p_n$ is uniquely recovered from the points $(p_i)_{i=1}^{n-1}$ and $o$. Since the latter variety is rational so long as $n\le 9$, it follows that $\Ycal_5(\mu)$ is rational when $\mu_n=1$, $n\le 9$ and also that $\dim \Ycal_5(\mu)=\dim \Zcal_5(\mu)+3=12+n$, as required.  
\vspace{0.5em}

If $\mu_n\ge 2$, then the possibilities are
$$\mu\in \{(2,2,3,3),(2,3,5),(3,3,4),(3,7),(5,5)\}$$
and the rationality question on $\Zcal_5(\mu)$ needs to be addressed case by case. The idea is simple: we turn the given relation on points into two relations, to which we apply the same method as in the case $\mu_n=1$. 
\vspace{0.5em}
\begin{enumerate}[label=$\roman*)$]
    \item If $\mu=(2,2,3,3)$, write $p_1+p_2+p_3+p_4\in |\Ocal_E(1)(x)|$ for a unique point $x\in E$, so the relation $\sum\limits_{i=1}^n\mu_ip_i+2o\in |\Ocal_E(4)|$ becomes $p_3+p_4+2o+2x\in |\Ocal_E(2)|$. These $2$ equations show that $p_1,p_4$ are uniquely recovered from $(E,p_2,p_3,o,x)$, and thus there is a birational isomorphism 
$$\Zcal_5(\mu)\dashrightarrow \{(E,p_2,p_3,o,x):E\in|\Ocal_{\PP^2}(3)|,p_2,p_3,o,x\in E\}$$
that proves the rationality of $\Zcal_5(\mu)$. 
    \item If $\mu=(2,3,5)$, one can define $x$ via $p_1+p_2+2p_3\in |\Ocal_E(1)(x)|$, in which case the relation becomes $p_2+p_3+2o+2x\in |\Ocal_E(2)|$, and thus there is a birational isomorphism 
$$\Zcal_5(\mu)\dashrightarrow \{(E,p_3,o,x):E\in |\Ocal_{\PP^2}(3)|\textnormal{, }p_3,o,x\in E\}$$
that, again, proves that $\Zcal_5(\mu)$ is rational. 
    \item If $\mu=(3,3,4)$, define $x$ via $p_1+p_2+p_3+o\in |\Ocal_E(x)|$ and this way there is a birational isomorphism 
$$\Zcal_5(\mu)\dashrightarrow \{(E,p_2,p_3,x):E\in |\Ocal_{\PP^2}(3)|\textnormal{, }p_2,p_3,x\in E\},$$
whence $\Zcal_5(\mu)$ is rational.
    \item If $\mu=(3,7)$, the point $x$ is defined by $p_1+2p_2+x\in |\Ocal_E(1)(x)|$, and then one obtains a birational isomorphism 
$$\Zcal_5(\mu)\dashrightarrow \{(E,p_2,x):E\in |\Ocal_{\PP^2}(3)|\textnormal{, }p_2,x\in E\}$$
that shows the rationality of $\Zcal_5(\mu)$. 
    \item If $\mu=(5,5)$, then define $x$ by $2p_1+p_2+o\in |\Ocal_E(1)(x)|$ and derive the birational isomorphism 
$$\Zcal_5(\mu)\dashrightarrow \{(E,p_2,x):E\in |\Ocal_{\PP^2}(3)|\textnormal{, }p_2,x\in E\}$$
that, again, shows $\Zcal_5(\mu)$ is rational.
\end{enumerate}

\textbf{(b)} \textit{Case of $\mu$ even and $\nu=(2)$.} Here, $\PP\Omega_5(\mu,\nu)$ consists of the two components $\PP\Omega_5(\mu,\nu)^+$ and $\PP\Omega_5(\mu,\nu)^-$. Let $\Xcal_5(\mu,\nu)^+$ and $\Xcal_5(\mu,\nu)^-$ be the corresponding loci in $\Xcal_5(\mu,\nu)$, respectively. These components are reflected in the geometry of $\Ycal_5(\mu)$, as we now explain. 
\vspace{0.5em}

Let $\left(E,o,(p_i)_{i=1}^n,\qbf\right)\in \Ycal_5(\mu)$ be a point, so $\sum\limits_{i=1}^n\mu_ip_i+2\qbf\in |\Ocal_E(6)|$. Since all $\mu_i$ are even, this implies that there exists a line bundle $\tau\in \Pic^0(E)$ such that $\tau^2=\Ocal_E$ and 
$$\sum_{i=1}^n\mu_i/2\cdot p_i+\qbf\in |\Ocal_E(3)\otimes \tau|.$$
If $\tau=\Ocal_E$, it follows that there exists a second cubic $E'\in |\Ocal_{\PP^2}(3)|$, $E'\ne E$, such that $E\cdot E'=\sum\limits_{i=1}^n\mu_i/2\cdot p_i+\qbf$. In these conditions, we obtain $\Gamma\cdot E'\ge \sum\limits_{i=1}^n\mu_i/2\cdot p_i+2\qbf$, so the remaining length $5$ scheme in their intersection needs to be linearly equivalent to $\sum\limits_{i=1}^n\mu_i/2\cdot p_i$ on $\Gamma$. This can only happen when $h^0\left(C,\sum\limits_{i=1}^n\mu_i/2\cdot p_i\right)\ge 2$, so this situation corresponds to the component $\PP\Omega_5(\mu,\nu)^-$. Defining
$$\Ycal_5(\mu)^-:=\left\{(E,(p_i)_{i=1}^n,\qbf):\begin{matrix}
    E\in |\Ocal_{\PP^2}(3)|\textnormal{ and }p_i\in E\textnormal{ with}\\
    \sum\limits_{i=1}^n\mu_i/2\cdot p_i+\qbf\in |\Ocal_E(3)|
\end{matrix}\right\}$$
and 
$$\Ycal_5(\mu)^+:=\left\{(E,\tau,(p_i)_{i=1}^n,\qbf):\begin{matrix}
    E\in |\Ocal_{\PP^2}(3)|\textnormal{ and }p_i\in E\textnormal{ with}\\
    \sum\limits_{i=1}^n\mu_i/2\cdot p_i+\qbf\in |\Ocal_E(3)\otimes \tau|\\
    \textnormal{for }\tau^2=\Ocal_E\textnormal{ with }\tau\ne \Ocal_E
\end{matrix}\right\},$$
these are exactly the connected components of $\Ycal_5(\mu)$ and, as we discussed earlier, there are $\PP^1$--bundles $\pi^-:\Xcal_5(\mu,-\nu)^-\to \Ycal_5(\mu)^-$ and $\pi^+:\Xcal_5(\mu,-\nu)^+\to \Ycal_5(\mu)^+$ over a dense open in the target. It is known that the varieties 
$$\Zcal_5(\mu)^-:=\left\{(E,(p_i)_{i=1}^n):E\in |\Ocal_{\PP^2}(3)|\textnormal{ and }p_i\in E\right\}$$
and 
$$\Zcal_5(\mu)^+:=\left\{(E,\tau,(p_i)_{i=1}^n):E\in |\Ocal_{\PP^2}(3)|,\textnormal{ }p_i\in E\textnormal{ and }\tau^2=\Ocal_E\textnormal{ for }\tau\ne \Ocal_E\right\}$$
are both rational, as $n\le 5$ (see \cite{krugPHD} and \cite{bellorouskyPHD}). Since the datum $\qbf$ is recovered via the relation $\qbf\in \left\lvert\Ocal_E(3)\otimes \tau\left(-\sum\limits_{i=1}^n\mu_i/2\cdot p_i\right)\right\rvert\cong \PP^{3}$, we conclude that both projections $\Ycal_5(\mu)^-\to \Zcal_5(\mu)^-$ and $\Ycal_5(\mu)^+\to \Zcal_5(\mu)^+$ realise the source as a dense open in a projective $\PP^3$--bundle over a dense open in the target. Consequently, both $\Ycal_5(\mu)^-$ and $\Ycal_5(\mu)^+$ are rational varieties, and also have dimension $12+n$. 
\vspace{0.5em}

\textbf{(c)} \textit{Case of $\mu$ even and $\nu=(1,1)$.} Here, $\PP\Omega_5(\mu,\nu)$ consists of two components $\PP\Omega_5(\mu,-\nu)^+$ and $\PP\Omega_5(\mu,\nu)^-$ as well, depending on the parity of $h^0\left(C,\sum\limits_{i=1}^n\mu_i/2\cdot p_i\right)$. Define $\Xcal_5(\mu,\nu)^+$ and $\Xcal_5(\mu,-\nu)^-$ as the corresponding incidence varieties. We can argue precisely as in \textbf{(b)}, to see that instead we obtain $\PP^1$--bundles $\pi^-:\Xcal_5(\mu,\nu)^-\to \Ycal_5(\mu)^+$ and $\pi^+:\Xcal_5(\mu,\nu)^+\to \Ycal_5(\mu)^-$ over a dense open in the target, where the varieties $\Ycal_5(\mu)^{\pm}$ are exactly those from the previous case. Consequently, $\Xcal_5(\mu,-\nu)^+$ and $\Xcal_5(\mu,-\nu)^-$ are again rational varieties. 
\end{proof}

\section{Meromorphic Differentials in Genus \texorpdfstring{$4$}{}}

We continue with the meromorphic differential strata in genus $4$. Theorem \ref{thm: unirat diff strata in genus 4} states that each component of the variety
$$\PP\Omega_4(\mu,\nu)=\left\{\left[C,(p_i)_{i=1}^n,(q_j)_{j=1}^r\right]\in\Mcal_{4,n+r}:\Ocal_C\left(\sum_{i=1}^n\mu_ip_i-\sum_{j=1}^r\nu_jq_j\right)\cong \omega_C\right\}$$
is unirational, whenever $s:=|\nu|\le 3$ and $n\le 8$. In what follows, the cases $s=2,3$ are treated together for most of the argument. 
\vspace{0.5em}

The canonical embedding of a nonhyperelliptic curve $C$ of genus $4$ presents it as the complete intersection of a quadric and a cubic in $\PP^3$. If we assume $C$ has no vanishing theta nulls, then the quadric is smooth. Set $S:=\PP^1\times \PP^1$ so that $C\in |\Ocal_S(3,3)|$. 

\begin{definition}
Consider the incidence variety 
$$\Xcal_4(\mu,\nu)=\left\{(C,E,Q,(p_i)_{i=1}^n,(q_j)_{j=1}^r,\qbf):\begin{matrix}
    C\in |\Ocal_S(3,3)|, E\in |\Ocal_S(2,2)|,\\
    Q\in |\Ocal_S(1,1)|\textnormal{ and }\qbf\in C_{6-s},\textnormal{ with}\\
    C\cdot E=\sum\limits_{i=1}^n\mu_ip_i+\qbf,C\cdot Q=\sum\limits_{j=1}^r\nu_jq_j+\qbf
\end{matrix}\right\}.$$
\end{definition}

As in the case of genus $5$, we have two similar lemmas that directly imply the unirationality of $\PP\Omega_4(\mu,\nu)$. 

\begin{lemma}\label{dominance lemma gen 4}
The map $\Xcal_4(\mu,\nu)\dashrightarrow \PP\Omega_4(\mu,\nu)$ defined by
$$\left[C,E,Q,(p_i)_{i=1}^n,(q_j)_{j=1}^r,\qbf\right]\mapsto \left[C,(p_i)_{i=1}^n,(q_j)_{j=1}^r\right]$$
is dominant over each component of the target, whenever $s\le 3$. Furthermore, this induces a $\PP^{3-s}$--bundle $\Xcal_4(\mu,\nu)\gitquot\Aut(S)\dashrightarrow \PP\Omega_4(\mu,\nu)$ that dominates each component of the target, whenever $s\le 3$. 
\end{lemma}
\begin{proof}
To see that the map is well-defined, notice that the intersection conditions of $C$ with $E$ and $Q$ give 
$$\sum_{i=1}^n\mu_ip_i+\qbf\in |\Ocal_C(2,2)|\textnormal{ and }\sum_{j=1}^r\nu_jq_j+\qbf\in |\Ocal_C(1,1)|, \textnormal{ so }\sum_{i=1}^n\mu_ip_i-\sum_{j=1}^r\nu_jq_j\sim K_C,$$
as required. Let $\left[C,(p_i)_{i=1}^n,(q_j)_{j=1}^r\right]\in \PP\Omega_4(\mu,\nu)$ be a point such that $C$ has no vanishing theta nulls (this can be achieved with models in $S$ described as in the definition of $\Xcal_4(\mu, \nu)$). Furthermore, we require a generality assumption on the points $q_j$, $1\le j\le r$, which we can also attain. In this case, there exists a linear system of conics $Q\in |\Ocal_S(1,1)|$, of dimension $3 - s$, such that $C\cdot Q\ge \sum\limits_{j=1}^r\nu_jq_j$. Let us write $C\cdot Q=\sum\limits_{j=1}^r\nu_jq_j+\qbf$, for some degree $6-s$ divisor $\qbf$ on $C$. By adjunction, we have $\omega_C=\Ocal_C(1,1)$, so the relation $\sum\limits_{i=1}^n\mu_ip_i-\sum\limits_{j=1}^r\nu_jq_j\sim K_C$ gives 
$$\sum_{i=1}^n\mu_ip_i+\qbf\in |\Ocal_C(2,2)|.$$
Since the restriction map $\H^0(S,\Ocal_S(2,2))\to \H^0(C,\Ocal_C(2,2))$ is an isomorphism, it follows that there exists a unique curve $E\in |\Ocal_S(2,2)|$ such that $C\cdot E=\sum\limits_{i=1}^n\mu_ip_i+\qbf$. The dominance statement is now proven. 
\end{proof}

The next lemma addresses the rationality of the incidence variety $\Xcal_4(\mu,\nu)$.

\begin{lemma}
All components of $\Xcal_4(\mu,\nu)$ are rational varieties, whenever $s\le 3$ and $n\le 8$. 
\end{lemma}
\begin{proof}
Take a point $\left[C,E,Q,(p_i)_{i=1}^n,(q_j)_{j=1}^r,\qbf\right]\in \Xcal_4(\mu,\nu)$ such that $p_g(C) = 4$. It follows that $E\cdot Q\ge \qbf$, so there exists $\rbf\in E_{s-2}$ such that $E\cdot Q=\qbf+\rbf$. Concretely, for $s=3$, $\rbf$ is a single point, while for $s=2$ it is empty. Introduce the variety 
\vspace{0.5em}
$$\Ycal_4(\mu,\nu):=\left\{\left(E,Q,(p_i)_{i=1}^n,(q_j)_{j=1}^r,\rbf\right):\begin{matrix}
    E\in |\Ocal_S(2,2)|, Q\in |\Ocal_S(1,1)|,\\
    E\cap Q\ge \rbf \textnormal{ for }\rbf\in E_{s-2},\,q_j\in Q,\\
    \sum\limits_{i=1}^n\mu_ip_i\in \left\lvert\Ocal_E(2,2)\left(\rbf\right)\right\rvert
\end{matrix}\right\},$$
so that there exists a map $\pi:\Xcal_4(\mu,\nu)\dashrightarrow \Ycal_4(\mu,\nu)$ that forgets the curve $C$. This is well-defined, since $\rbf$ and $\dbf$ determine each other uniquely. As was the case for genus $5$, we split the discussion into the following three claims, since they together ensure the rationality of $\Xcal_4(\mu,\nu)$:
\begin{enumerate}[label=\textbf{(\roman*)}]
    \item\label{indep lemma gen 4} the map $\pi$ is a $\PP^{4-s}$--bundle over its image in any component of the target,
    \item\label{rat check for Ycal4} all components of $\Ycal_4(\mu,\nu)$ are rational if $n\le 8$ and $s\le 3$, and 
    \item\label{dim count gen 4} $\dim \Xcal_4(\mu,\nu)=14+n+r-s$ and $\dim \Ycal_4(\mu)=10+n+r$, so $\pi$ is, in fact, a $\PP^{4-s}$--bundle over a dense open in each component of the target. 
\end{enumerate}

Theorem \ref{big indep lemma} implies part \ref{indep lemma gen 4}. Indeed, the divisor of conditions on $|\Ocal_S(3,3)|$ is $Z=\sum\limits_{i=1}^n\mu_ip_i+\sum\limits_{j=1}^r\nu_jq_j+\qbf$, and on the curve $C$ we have $Z=2K_C+\sum\limits_{j=1}^r\nu_jq_j$, so the number of independent conditions imposed on the system $|\Ocal_S(3,3)|\cong \PP^{15}$ is
\[
\length(Z)-h^0\left(C,\sum\limits_{j=1}^r\nu_jq_j\right)=(12+s)-1=11+s.
\]
\vspace{0.5em}

Furthermore, Lemma \ref{dominance lemma gen 4} shows in particular that $\dim \Xcal_4(\mu,\nu)=14+n+r-s$. The second dimensionality check in \ref{dim count gen 4} is treated together with \ref{rat check for Ycal4} in the following casework. 

\vskip 0.5em
\textbf{(a)} \textit{Case of odd $\mu$, when $s=2$.} Here, $\PP\Omega_4(\mu,\nu)$ consists of a single component. The same holds for $\Ycal_4(\mu,\nu)$, as it is transparent from the argument below. We have 
\begin{align*}
    \Ycal_4(\mu,\nu)&=\left\{\left(E,(p_i)_{i=1}^n\right):E\in |\Ocal_S(2,2)|,\sum\limits_{i=1}^n\mu_ip_i\in |\Ocal_E(2,2)|\right\}\times\\
    &\times\left\{\left(Q,(q_j)_{j=1}^r\right):Q\in |\Ocal_S(1,1)|, q_j\in Q\right\}
\end{align*}
and both factors are rational. This is immediate for the latter, while the former needs to be separated into further subcases.
\vspace{0.5em}

\begin{enumerate}[label=$\roman*)$]
    \item If $\mu_n=1$, then as soon as the points $(p_i)_{i=1}^{n-1}\in S^{n-1}$ determine an elliptic curve $E\in |\Ocal_S(2,2)|$ passing through them, $p_n$ is uniquely determined by the condition $\sum\limits_{i=1}^n\mu_ip_i\in |\Ocal_E(2,2)|$.
    \item If $\mu=(2,3,3)$, the same method as in genus $5$ leads to success. There is a unique point $x\in E$ such that $p_1+p_2+p_3\in |\Ocal_E(1,0)(x)|$, and hence the condition $2p_1+3p_2+3p_3\in |\Ocal_E(2,2)|$ is equivalent to $p_2+p_3+2x\in |\Ocal_E(0,2)|$. Therefore, the first factor of $\Ycal_4(\mu,\nu)$ is birationally isomorphic to 
$$\{(E,p_3,x):E\in |\Ocal_S(2,2)|\textnormal{, }p_3,x\in E\}$$
which is obviously rational.
    \item If $\mu=(3,5)$, define the point $x\in E$ such that $p_1+2p_2\in |\Ocal_E(1,0)(x)|$, and notice that the condition $3p_1+5p_2\in |\Ocal_E(2,2)|$ is equivalent to $3x\in |\Ocal_E(-1,2)(p_2)|$, so the first factor of $\Ycal_4(\mu,\nu)$ is birationally isomorphic to 
$$
\{(E,x):E\in |\Ocal_S(2,2)|\textnormal{, }x\in E\}\textnormal{, clearly a rational variety.}
$$
\end{enumerate}

\textbf{(b)} \textit{Case of $\mu$ even, when $\nu=(2)$.} The variety $\Ycal_4(\mu,\nu)$ has the same expression as in \textbf{(a)}. Define $\Xcal_4(\mu,\nu)^+$ and $\Xcal_4(\mu,\nu)^-$ to be the loci corresponding to the two components $\PP\Omega_4(\mu, \nu)^+$ and $\PP\Omega_4(\mu, \nu)^-$ of $\PP\Omega_4(\mu, \nu)$. We expect that $\Ycal_4(\mu)$ also has $2$ components and we check this below. 
\vspace{0.5em}

Let $\Zcal_4(\mu)$ denote the first factor of $\Ycal_4(\mu,\nu)$. Consider a point $(E,(p_i)_{i=1}^n)\in \Zcal_4(\mu)$. Since $\mu_i$ are all even, the condition $\sum\limits_{i=1}^n\mu_ip_i\in |\Ocal_E(2,2)|$ implies that there exists a line bundle $\tau$ on $E$ such that $\tau^2=\Ocal_E$ and 
$$\sum_{i=1}^n\mu_i/2\cdot p_i\in |\Ocal_E(1,1)\otimes \tau|.$$
If $\tau=\Ocal_E$, then there exists $Q'\in |\Ocal_S(1,1)|$ such that $E\cdot Q'=\sum\limits_{i=1}^n\mu_i/2\cdot p_i$. In particular, $C\cdot Q'\ge \sum\limits_{i=1}^n\mu_i/2\cdot p_i$, so the remaining length $2$ scheme in their intersection is linearly equivalent to $\sum\limits_{i=1}^n\mu_i/2\cdot p_i-q$, meaning $h^0\left(C,\sum\limits_{i=1}^n\mu_i/2\cdot p_i-q\right)\ge 1$ is necessary. Therefore, this case corresponds to the component $\PP\Omega_4(\mu,\nu)^-$. We deduce that 
$$\Zcal_4(\mu)^-:=\left\{(E,(p_i)_{i=1}^n):E\in |\Ocal_S(2,2)|\textnormal{ and }\sum\limits_{i=1}^n\mu_i/2\cdot p_i\in |\Ocal_E(1,1)|\right\}$$
and 
$$\Zcal_4(\mu)^+:=\left\{(E,\tau,(p_i)_{i=1}^n):\begin{matrix}
    E\in |\Ocal_S(2,2)|\textnormal{ and }\sum\limits_{i=1}^n\mu_i/2\cdot p_i\in |\Ocal_E(1,1)\otimes \tau|\\
    \textnormal{for some }\tau\ne \Ocal_E\textnormal{ with }\tau^2=\Ocal_E
\end{matrix}\right\}$$
are the irreducible components of $\Zcal_4(\mu)$. These in turn define the irreducible components of $\Ycal_4(\mu,\nu)$, which we shall call $\Ycal_4(\mu,\nu)^-$ and $\Ycal_4(\mu,\nu)^+$ accordingly. We have maps 
$$\pi^-:\Xcal_4(\mu,\nu)^-\to \Ycal_4(\mu,\nu)^-\textnormal{ and }\pi^+:\Xcal_4(\mu,\nu)^+\to \Ycal_4(\mu,\nu)^+$$
that are $\PP^2$-bundles over a dense open in the target. The varieties $\Zcal_4(\mu)^-$ and $\Zcal_4(\mu)^+$ are rational, as the following subcases show: 

\begin{enumerate}[label=$\roman*)$]
    \item If $\mu_n=2$, then one can use the fact that $\Rcal_{1,n}$ is rational for $n\le 6$ (see \cite{krugPHD}), and then argue just as in the $\mu_n=1$ part of case \textbf{(a)}. 
    \item If $\mu=(4,4)$, then we address the rationality question of $\Zcal_4(\mu)^+$ by applying the same arguments as we proposed for the unirationality question of the spaces $\PP\Omega_g(\mu,\nu)$. Concretely, one defines the incidence variety 
$$\Tcal^+:=\left\{\left(E,(Q_i)_{i=1}^3,p_1,p_2,x_1,x_2\right):\begin{matrix}
    E\in |\Ocal_S(2,2)|\textnormal{, }Q_i\in |\Ocal_S(1,1)|\\
    E\cdot Q_1=p_1+2p_2+x_1, E\cdot Q_2=3p_1+x_2\\
    \textnormal{and }E\cdot Q_3=p_1+2x_1+x_2
\end{matrix}\right\},$$
which is rational and admits a birational isomorphism $\Tcal^+\dashrightarrow \Zcal_4(\mu)^+$ defined by $$\left[E,(Q_i)_{i=1}^3,p_1,p_2,x_1,x_2\right]\mapsto \left[E,\Ocal_E(x_1-p_1),(p_i)_{i=1}^2\right].$$ We omit the details, and the proof for $\Zcal_4(\mu)^-$ is treated analogously. The case $\mu=(8)$ can be dealt with in a similar manner. 
\end{enumerate}

\textbf{(c)} \textit{Case of $\mu$ even, when $\nu=(1,1)$.} As in the case of genus $5$, the change from $\nu=(2)$ to $\nu=(1,1)$ is immediate and we omit the details. 
\vspace{0.5em}

\textbf{(d)} \textit{Case of $n\le 8$, when $s=3$.} The variety $\PP\Omega_4(\mu,\nu)$ is irreducible in this case. We have 
$$\Ycal_4(\mu,\nu):=\left\{(E,Q,(p_i)_{i=1}^n,(q_j)_{j=1}^r,\rbf):\begin{matrix}
    E\in |\Ocal_S(2,2)|, Q\in |\Ocal_S(1,1)|, Q\cdot E\ge \rbf\\
    \textnormal{with }q_j\in Q\textnormal{ and }\sum\limits_{i=1}^n\mu_ip_i\in |\Ocal_E(2,2)(\rbf)|
\end{matrix}\right\}.$$
Notice that a general tuple $(p_i)_{i=1}^n\in S^n$ gives a $\PP^{8-n}$--bundle of elliptic curves $E\in |\Ocal_S(2,2)|$ passing through the $p_i$. The point $\rbf$ is uniquely determined from the condition $\sum\limits_{i=1}^n\mu_ip_i\in |\Ocal_E(2,2)(\rbf)|$, and the data of a conic $Q\in |\Ocal_S(1,1)|$ passing through $\rbf$ and having $r$ marked points $(q_j)_{j=1}^r$ is parametrized by $\textnormal{Stab}_{\PGL(2)}(\rbf)\times \left(\PP^1\right)^r\cong \GG_m\times \AA^1\times \left(\PP^1\right)^r$. These show that $\Ycal_4(\mu,\nu)$ is a rational variety. 
\end{proof}

\section{Meromorphic differential strata in genus \texorpdfstring{$3$}{}}

In this section  we study the meromorphic differential strata in genus $3$. Theorem \ref{thm: unirat diff strata in genus 3} states that the components of the variety
$$
\PP\Omega_3(\mu,\nu) = \left\{\left[C,(p_i)_{i=1}^n,(q_j)_{j=1}^r\right] \in \Mcal_{3,n +r}\colon \Ocal_C\left(\sum_{i=1}^n \mu_i p_i - \sum_{j=1}^r \nu_jq_j\right) \cong \omega_C\right\}
$$
are unirational, when the number of poles $s:=|\nu|$ is at most $4$. 

The method requires that we study each of the following cases individually: $s=2,4$ and $\nu\in \{(1,1),(2),(2,2),(4)\}$, $s=3,4$ and $\nu\in \{(1,1,1),(1,2),(3),(1,1,1,1),(1,1,2),(1,3)\}$, i.e. all the other possible values. 

\subsection{Meromorphic strata with at most $2$ poles with equal orders}

This subsection is dedicated to the analysis of the cases when 
\begin{align}\label{values of nu when sum nu 2 or even 4}
    \nu\in \{(1,1),(2),(2,2),(4)\}.
\end{align}
Let us write $\nu=\varepsilon \cdot \nu'$, where $\varepsilon\in \{1,2\}$ and $\nu'\in \{(1,1),(2)\}$, and let $r$ be the number of distinct poles, so $r\in \{2,1\}$. The spaces of interest are 
$$\PP\Omega_3(\mu,\nu)=\left\{\left[C,(p_i)_{i=1}^n,q^{\pm}\right]\in \Mcal_{3,n+r}:\Ocal_C\left(\sum_{i=1}^n\mu_ip_i-\varepsilon\left(q^{+}+q^-\right)\right)\cong \omega_C\right\}$$
where, if $\nu\in\{(2),(4)\}$, the notation stands for $q^+=q^-=:q$. If $C$ is a nonhyperelliptic curve of genus $3$, the canonical embedding presents it as a plane quartic $C\in |\Ocal_{\PP^2}(4)|$. 

\begin{remark}
Consider for the moment the case when $\mu$ is even and $\nu\in \{(1,1),(2)\}$. In this case, $\PP\Omega_3(\mu,\nu)$ consists of the irreducible components $\PP\Omega_3(\mu,\nu)^+$ and $\PP\Omega_3(\mu,\nu)^-$. Observe that, since the point $q$ is distinct from the $p_i$ and $\sum\limits_{i=1}^n\mu_i/2\cdot p_i$ is an effective degree $3$ divisor, the description of the irreducible components in both cases with $|\nu|=2$ reduces to whether $\sum\limits_{i=1}^n\mu_i/2\cdot p_i$ moves or does not move. 
\end{remark}

We split the analysis into the following cases: 

\begin{enumerate}[label=\textbf{(\alph*)}]
    \item\label{sum v 2, m even and moving div} \textit{Case of $|\nu|=2$, $\mu$ even and moving $\sum\limits_{i=1}^n\mu_i/2\cdot p_i$}, corresponding to the components $\PP\Omega_3(\mu,\nu)^-$ for $\nu=(2)$ and $\PP\Omega_3(\mu,\nu)^+$ for $\nu=(1,1)$, 
    \item\label{sum v 2, arbitrary m} \textit{Case of $|\nu|=2$ and any other $\mu$}, which for $\mu$ even corresponds to the component $\PP\Omega_3(\mu,\nu)^+$ for $\nu=(2)$ and $\PP\Omega_3(\mu,\nu)^-$ for $\nu=(1,1)$, and for any other $\mu$ has $\PP\Omega_3(\mu,\nu)$ irreducible,
    \item\label{sum v 4 & even, mu_n=1} \textit{Case of $\nu\in \{(2,2),(4)\}$ and odd $\mu$}, corresponding to an irreducible $\PP\Omega_3(\mu,\nu)$, and
    \item\label{sum v 4 & even, m even} \textit{Case of $\nu\in \{(2,2),(4)\}$ and $\mu$ even}, corresponding to two irreducible components $\PP\Omega_3(\mu,\nu)^+$ and $\PP\Omega_3(\mu,\nu)^-$. 
\end{enumerate}
\vspace{0.5em}
Situation \ref{sum v 2, m even and moving div} can be handled from the very beginning. The remaining other $3$ situations are treated together afterwards, until their splitting in Lemma \ref{rational variety w nu even & sum 2 or 4}. 
\vspace{0.5em}

\ref{sum v 2, m even and moving div} \textit{Case of $|\nu|=2$, $\mu$ even and moving $\sum\limits_{i=1}^n\mu_i/2\cdot p_i$.} This corresponds to $\PP\Omega_3(\mu,\nu)^-$ for $\nu=(2)$ and to $\PP\Omega_3(\mu,\nu)^+$ for $\nu=(1,1)$. We need the following 

\begin{definition}\label{incidence var for nu 2 and 1,1 w moving div}
Consider the incidence variety 
$$\Xcal_3(\mu,\nu)=\left\{(C,L,\ell,(p_i)_{i=1}^n,q^{\pm},x):\begin{matrix}
    C\in |\Ocal_{\PP^2}(4)|\textnormal{ and } L,\ell\in \left(\PP^2\right)^\vee\textnormal{ with }\\
    C\cdot L=\sum\limits_{i=1}^n\mu_i/2\cdot p_i+x, C\cdot \ell=q^++q^-+2x
\end{matrix}\right\}.$$
When $\nu=(2)$, the notation stands for $q^{+}=q^-=:q$. 
\end{definition}

As reasoned previously, the following two lemmas are enough to ensure the unirationality of $\PP\Omega_3(\mu,\nu)^-$ for $\nu=(2)$ and of $\PP\Omega_3(\mu,\nu)^+$ for $\nu=(1,1)$. 

\begin{lemma}\label{dominant map for nu 2 and 1,1 w moving div}
The maps $\Xcal_3(\mu,\nu)\dashrightarrow \PP\Omega_3(\mu,\nu)^-$ for $\nu=(2)$ and $\Xcal_3(\mu,\nu)\dashrightarrow\PP\Omega_3(\mu,\nu)^+$ for $\nu=(1,1)$ defined by 
$$\left[C,L,\ell,(p_i)_{i=1}^n,q^{\pm},x)\right]\mapsto \left[C,(p_i)_{i=1}^n,q^{\pm}\right]$$
are well-defined and dominant. 
\end{lemma}
\begin{proof}
We prove this only when $\nu=(1,1)$, with the case $\nu=(2)$ being completely similar. Pick a point $\left[C,(p_i)_{i=1}^n,q^{\pm}\right]\in \PP\Omega_3(\mu,\nu)^+$ such that $C$ is nonhyperelliptic and such that $\sum\limits_{i=1}^n\mu_i/2\cdot p_i$ is a moving divisor. Then there exists a point $x\in C$ and a line $L$ such that 
$$C\cdot L=\sum_{i=1}^n\mu_i/2\cdot p_i+x.$$
Using the relation $\sum\limits_{i=1}^n\mu_ip_i-q^+-q^-\sim K_C$, we deduce that there also exists a line $\ell$ such that $C\cdot \ell=q^++q^-+2x$. This concludes the proof of the lemma. 
\end{proof}

\begin{lemma}\label{rat variety for nu 2 and 1,1 w moving div}
The variety $\Xcal_3(\mu,\nu)$ is rational.
\end{lemma}
\begin{proof}
Introduce the variety 
$$\Ycal_3(\mu,\nu):=\left\{(L,(p_i)_{i=1}^n):L\textnormal{ line, }p_i\in L\right\}\times \{(\ell,q^{\pm}):\ell\textnormal{ line, }q^{\pm}\in \ell\}$$
which is clearly rational and notice that the forgetful map $\pi:\Xcal_3(\mu,\nu)\to \Ycal_3(\mu,\nu)$ is a $\PP^{7}$--bundle over the image. Indeed, the divisor of conditions on $|\Ocal_{\PP^2}(4)|$ is given by $Z=\sum\limits_{i=1}^n\mu_i/2\cdot p_i+q^++q^-+2x$, which on $C$ is $Z=2K_C-x$, so there are $\length(Z)-h^0(C,x-K_C)=7$ independent linear conditions imposed on the system $|\Ocal_{\PP^2}(4)|\cong \PP^{14}$, according to Theorem \ref{big indep lemma}.

Finally, since $\dim \Xcal_3(\mu,\nu)=11+n+r$ and $\dim \Ycal_3(\mu,\nu)=4+n+r$, it follows that $\pi$ is a $\PP^7$--bundle over a dense open in the target. This concludes the proof of the lemma. 
\end{proof}

Putting together Lemma \ref{dominant map for nu 2 and 1,1 w moving div} and Lemma \ref{rat variety for nu 2 and 1,1 w moving div}, we infer that the components $\PP\Omega_3(\mu,\nu)^-$ for $\nu=(2)$ and $\PP\Omega_3(\mu,\nu)^+$ for $\nu=(1,1)$ are unirational. This settles the case \ref{sum v 2, m even and moving div} and we may proceed with the collective analysis of the cases \ref{sum v 2, arbitrary m}, \ref{sum v 4 & even, mu_n=1} and \ref{sum v 4 & even, m even}. 

\begin{definition}
For $\mu$ and $\nu$ as in \ref{sum v 2, arbitrary m}, \ref{sum v 4 & even, mu_n=1} or \ref{sum v 4 & even, m even}, consider the incidence variety 
$$\Xcal_3(\mu,\nu)=\left\{\left(C,E,\ell,(p_i)_{i=1}^n,q^{\pm},\qbf\right):\begin{matrix}
    C\in |\Ocal_{\PP^2}(4)|, E\in |\Ocal_{\PP^2}(1+\varepsilon)|, \ell\in \left(\PP^2\right)^\vee,\\
    C\cdot E=\sum\limits_{i=1}^n\mu_ip_i+\varepsilon\qbf\textnormal{ and }C\cdot \ell=q^++q^-+\qbf\\
    \textnormal{for some }\qbf\in C_2
\end{matrix}\right\}.$$
In the cases $\nu\in\{(2),(4)\}$, the notation $q^{\pm}$ stands for $q^+=q^-=:q$. 
\end{definition}

As in the previous sections, there are two lemmas that imply the unirationality of all components of $\PP\Omega_3(\mu,\nu)$. 

\begin{lemma}\label{dominance lemma gen 3 & 2 distinct poles}
Let $\mu$ and $\nu$ be as in \ref{sum v 2, arbitrary m}, \ref{sum v 4 & even, mu_n=1} or \ref{sum v 4 & even, m even}. The map $\Xcal_3(\mu,\nu)\dashrightarrow \PP\Omega_3(\mu,\nu)$ defined by 
$$\left[C,E,\ell,(p_i)_{i=1}^n,q^{\pm},\qbf\right]\mapsto \left[C,(p_i)_{i=1}^n,q^{\pm}\right]$$
is well-defined and dominates each component of the target. Furthermore, this induces a birational isomorphism $\Xcal_3(\mu,\nu)\gitquot \PGL(3)\dashrightarrow \PP\Omega_3(\mu,\nu)$.
\end{lemma}
\begin{proof}
To see that this map is well-defined, observe that the intersection conditions of $C$ with $E$ and $\ell$ give 
$$\sum_{i=1}^n\mu_ip_i+\varepsilon\qbf\in |\Ocal_C(1+\varepsilon)|\textnormal{ and }q^++q^-+\qbf\in |\Ocal_C(1)|,\textnormal{ so }\sum_{i=1}^n\mu_ip_i-\varepsilon(q^++q^-)\sim K_C.$$
Now we address the dominance claim. The existence of at most $2$ poles now enters the picture: one can cover the data of the poles with a single line. Indeed, there exists a unique line $\ell$ such that $C\cdot \ell\ge q^++q^-$, and one can write 
$$C\cdot \ell=q^++q^-+\qbf,\textnormal{ for some degree }2\textnormal{ divisor }\qbf.$$
The condition $\sum\limits_{i=1}^n\mu_ip_i-\varepsilon(q^++q^-)\sim K_C$ gives $\sum\limits_{i=1}^n\mu_ip_i+\varepsilon\qbf\in |\Ocal_C(\varepsilon+1)|$. Since the restriction map $\H^0(\PP^2,\Ocal_{\PP^2}(1+\varepsilon))\to \H^0(C,\Ocal_C(1+\varepsilon))$ is an isomorphism, it follows that there exists a unique curve $E\in |\Ocal_{\PP^2}(1+\varepsilon)|$ such that 
$C\cdot E=\sum\limits_{i=1}^n\mu_ip_i+\varepsilon\qbf$. This establishes the dominance claim. 
\end{proof}

As in the previous sections, the second lemma is concerned with the rationality of all components of $\Xcal_3(\mu,\nu)$.

\begin{lemma}\label{rational variety w nu even & sum 2 or 4}
All components of $\Xcal_3(\mu,\nu)$ are rational varieties when $\mu$ and $\nu$ are as in \ref{sum v 2, arbitrary m}, \ref{sum v 4 & even, mu_n=1} or \ref{sum v 4 & even, m even}. 
\end{lemma}
\begin{proof}
Consider the variety
$$\Ycal_3(\mu,\nu)=\left\{\left(E,\ell,(p_i)_{i=1}^n,q^{\pm},\qbf\right):\begin{matrix}
    E\in |\Ocal_{\PP^2}(1+\varepsilon)|, \ell\in \left(\PP^2\right)^\vee, \qbf\in E_2,q^{\pm}\in \ell,\\
    \sum\limits_{i=1}^n\mu_ip_i+\varepsilon\qbf\in |\Ocal_E(4)|\textnormal{ and }E\cdot \ell\ge \qbf
\end{matrix}\right\}.$$
There is a map $\pi:\Xcal_3(\mu,\nu)\to \Ycal_3(\mu,\nu)$ that forgets the curve $C$. We claim the following $3$ statements, which directly imply that all components of $\Xcal_3(\mu,\nu)$ are rational:
\begin{enumerate}[label=\textbf{(\roman*)}]
    \item\label{indep lemma gen 3 & 2 distinct poles} the map $\pi$ is a $\PP^{4}$--bundle over the image when $\varepsilon=1$ and a $\PP^1$--bundle over the image when $\varepsilon=2$, 
    \item\label{rat check for Ycal3 & 2 distinct poles} all components of the variety $\Ycal_3(\mu,\nu)$ are rational, when $\mu$ and $\nu$ are as in \ref{sum v 2, arbitrary m}, \ref{sum v 4 & even, mu_n=1} or \ref{sum v 4 & even, m even}, 
    \item\label{dim count gen 3 & 2 distinct poles} we have $\dim \Xcal_3(\mu,\nu)=11+n+r$ and $\dim \Ycal_3(\mu,\nu)$ is $7+n+r$ for $\varepsilon=1$ and $10+n+r$ for $\varepsilon=2$. 
\end{enumerate}

According to Theorem \ref{big indep lemma}, the divisor of conditions is $Z=\sum\limits_{i=1}^n\mu_ip_i+q^{+}+q^-+\varepsilon\qbf$, which on $C$ is $Z=(1+\varepsilon)K_C+q^++q^-$, so the number of independent linear conditions imposed on the system $|\Ocal_{\PP^2}(4)|\cong \PP^{14}$ is
\[
\length(Z)-h^0(C,Z-3K_C)=4\varepsilon+6-h^0(C,(\varepsilon-2)K_C+q^++q^-)\in \{10,13\},
\]
depending on whether $\varepsilon=1$ or $\varepsilon=2$. This establishes part \ref{indep lemma gen 3 & 2 distinct poles}.
\vspace{0.5em}

Lemma \ref{dominance lemma gen 3 & 2 distinct poles} shows that there exists a birational isomorphism $\Xcal_3(\mu,\nu)/\PGL(3)\dashrightarrow \PP\Omega_3(\mu,\nu)$, and thus $\dim \Xcal_3(\mu,\nu)=11+n+r$. The second part of \ref{dim count gen 3 & 2 distinct poles} and the statement \ref{rat check for Ycal3 & 2 distinct poles} are treated together, in the casework below. 
\vspace{0.5em}

\ref{sum v 2, arbitrary m} \textit{Case of $|\nu|=2$ and any other $\mu$.} Here, we have $\varepsilon=1$ and we can write 
$$\Ycal_3(\mu,\nu)=\left\{\left(E,(p_i)_{i=1}^n\right):E\textnormal{ conic, }p_i\in E\right\}\times \left\{\left(\ell,q^{\pm}\right):\ell\textnormal{ line, }q^{\pm}\in \ell\right\},$$
since we recover $\qbf$ via $\qbf=E\cdot \ell$ and condition $\sum\limits_{i=1}^n\mu_ip_i+\qbf\in |\Ocal_E(4)|$ is vacuous in this case. It is now clear that $\Ycal_3(\mu,\nu)$ is an irreducible rational variety.
\vspace{0.5em}

For the remaining cases \ref{sum v 4 & even, mu_n=1} and \ref{sum v 4 & even, m even}, we have $\varepsilon=2$, so $E$ is a plane cubic and $\length(E\cap \ell)=3$. It follows that there exists a point $x\in \PP^2$ such that $E\cdot \ell=\qbf+x$, so the condition $\sum\limits_{i=1}^n\mu_ip_i+2\qbf\in |\Ocal_E(4)|$ can be rewritten into 
$$\sum_{i=1}^n\mu_ip_i\in |\Ocal_E(2)(2x)|.$$

\ref{sum v 4 & even, mu_n=1} \textit{Case of $\nu\in\{(2,2),(4)\}$ and odd $\mu$.} From the above, we can write 
$$\Ycal_3(\mu,\nu)=\left\{\left(E,\ell,(p_i)_{i=1}^n,x,q^{\pm}\right):\begin{matrix}
    E\in |\Ocal_{\PP^2}(3)|, \ell\in \left(\PP^2\right)^\vee, x\in E\cap \ell,\\
    q^{\pm}\in \ell\textnormal{ and }\sum\limits_{i=1}^n\mu_ip_i\in |\Ocal_E(2)(2x)|
\end{matrix}\right\}.$$
The cases $\mu\in \{(2,3,3),(3,5)\}$ are treated exactly as the cases $\mu\in\{(2,2,3,3),(2,3,5)\}$ of genus $5$, and hence we omit them. This leaves us with the odd partitions that satisfy $\mu_n=1$. Since $n\le 8$, the variety $\Zcal_3(\mu,\nu):=\left\{(E,(p_i)_{i=1}^{n-1},x):E\textnormal{ cubic, }p_i,x\in E\right\}$ is rational. Also, the projection $\Ycal_3(\mu,\nu)\to \Zcal_3(\mu,\nu)$ that forgets the points $p_n$, $q^{\pm}$ and the line $\ell$ is a tower of projective bundles: the point $p_n$ is uniquely recovered via the condition $\sum\limits_{i=1}^n\mu_ip_i\in |\Ocal_E(2)(2x)|$, whereas $\ell$ varies in the pencil of lines through $x$, and $q^{\pm}$ are $r$ points on $\ell$. Consequently, $\Ycal_3(\mu,\nu)$ is an irreducible rational variety, and therefore so is $\Xcal_3(\mu,\nu)$. 

\vspace{0.5em}
\ref{sum v 4 & even, m even} \textit{Case of $\nu\in \{(2,2),(4)\}$ and $\mu$ even.} We recall that there are two components $\PP\Omega_3(\mu,\nu)^+$ and $\PP\Omega_3(\mu,\nu)^-$ in this case.

The variety $\Ycal_3(\mu,\nu)$ has the same expression as in \ref{sum v 4 & even, mu_n=1}. The condition $\sum\limits_{i=1}^n\mu_ip_i\in |\Ocal_E(2)(2x)|$ implies that there exists a $2$-torsion bundle $\tau$ on $E$ such that 
$$\sum_{i=1}^n\mu_i/2\cdot p_i\in |\Ocal_E(1)(x)\otimes \tau|.$$
If $\tau=\Ocal_E$, we recover $\qbf\in E_2$ via $E\cdot \ell=x+\qbf$ and we have $\qbf+\sum\limits_{i=1}^n\mu_i/2\cdot p_i\in |\Ocal_E(2)|$. This implies that there exists $Q'\in |\Ocal_{\PP^2}(2)|$ such that $E\cdot Q'=\qbf+\sum\limits_{i=1}^n\mu_i/2\cdot p_i$, which in turn gives $C\cdot Q'=\sum\limits_{i=1}^n\mu_i/2\cdot p_i+\qbf+\rbf$, where $\rbf\in C_2$ is some degree $2$ effective divisor. Comparing with the relation $\sum\limits_{i=1}^n\mu_ip_i+2\qbf\in |\Ocal_C(3)|$, we obtain $$2\rbf\in |\Ocal_C(1)|\textnormal{ and }\sum\limits_{i=1}^n\mu_i/2\cdot p_i-q^+-q^-\sim \rbf\textnormal{, so }h^0\left(\sum\limits_{i=1}^n\mu_i/2\cdot p_i-q^+-q^-\right)\ge 1.$$ 
We conclude that $\tau=\Ocal_E$ corresponds to the component $\PP\Omega_3(\mu,\nu)^-$, whereas $\tau\ne \Ocal_E$ correspond to the other component $\PP\Omega_3(\mu,\nu)^+$. These also define the two irreducible components of $\Ycal_3(\mu,\nu)$, namely 
$$\Ycal_3(\mu,\nu)^-=\left\{\left(E,\ell,(p_i)_{i=1}^n,x,q^{\pm}\right):\begin{matrix}
    E\in |\Ocal_{\PP^2}(3)|, \ell\in \left(\PP^2\right)^\vee, x\in E\cap \ell,\\
    q^{\pm}\in \ell\textnormal{ and }\sum\limits_{i=1}^n\mu_i/2\cdot p_i\in |\Ocal_E(1)(x)|
\end{matrix}\right\}$$
and
$$\Ycal_3(\mu,\nu)^+=\left\{\left(E,\tau,\ell,(p_i)_{i=1}^n,x,q^{\pm}\right):\begin{matrix}
    E\in |\Ocal_{\PP^2}(3)|, \ell\in \left(\PP^2\right)^\vee, x\in E\cap \ell,\\
    q^{\pm}\in \ell\textnormal{ and }\sum\limits_{i=1}^n\mu_i/2\cdot p_i\in |\Ocal_E(1)(x)\otimes \tau|\\
    \textnormal{for some }\tau\ne \Ocal_E\textnormal{ with }\tau^2\cong \Ocal_E
\end{matrix}\right\}.$$
Since in this case $n\le 4$, thanks to the fact that $\mu$ is even, one concludes that both these components are rational (see \cite{krugPHD}). These establish the rationality of the corresponding irreducible components $\Xcal_3(\mu,\nu)^-$ and $\Xcal_3(\mu,\nu)^+$ in $\Xcal_3(\mu,\nu)$, and finally grant the unirationality of all the components of $\PP\Omega_3(\mu,\nu)$ in this case as well. 
\end{proof}

\subsection{Meromorphic strata with $|\nu|=3,4$ and $\nu\ne (2,2),(4)$}

In this situation, all possible partitions of the poles are 
\begin{align}\label{values of nu for irred strata in sum 3,4}
    \nu\in \{(1,1,1),(1,2),(3),(1,1,1,1),(1,1,2),(1,3)\}.
\end{align}
For any of them, the meromorphic stratum of differentials $\PP\Omega_3(\mu,\nu)$ is irreducible. As in the previous subsection, we work with the canonical embedding for nonhyperelliptic curves $[C]\in \Mcal_3$.

\begin{definition}
Set $s=|\nu|\in\{3,4\}$ and consider the incidence variety 
$$\Xcal_3(\mu,\nu)=\left\{\left(C,E,Q,(p_i)_{i=1}^n,(q_j)_{j=1}^r,\qbf,\rbf\right):\begin{matrix}
    C\in |\Ocal_{\PP^2}(4)|, \textnormal{ }E\in |\Ocal_{\PP^2}(3)|,\textnormal{ }\qbf\in C_{8-s}\\
    \textnormal{with }C\cdot E=\sum\limits_{i=1}^n\mu_ip_i+\qbf,\textnormal{ }Q\in |\Ocal_{\PP^2}(2)|\\
    \textnormal{with }C\cdot Q=\sum\limits_{j=1}^r\nu_jq_j+\qbf,\textnormal{ }\rbf\in E_{s-2}\\
    \textnormal{with }E\cdot Q=\rbf+\qbf
\end{matrix}\right\}.$$
\end{definition}

We claim that $\Xcal_3(\mu,\nu)$ is a rational parametrisation for $\PP\Omega_3(\mu,\nu)$, whenever $\nu$ is as in (\ref{values of nu for irred strata in sum 3,4}). This is a consequence of the following two lemmas. 

\begin{lemma}\label{dominance in genus 3 for odd nu 3, 4}
The map $\Xcal_3(\mu,\nu)\dashrightarrow \PP\Omega_3(\mu,\nu)$ given by 
$$\left[C,E,Q,(p_i)_{i=1}^n,(q_j)_{j=1}^r,\qbf,\rbf\right]\mapsto \left[C,(p_i)_{i=1}^n,(q_j)_{j=1}^r\right]$$
is well-defined and dominant. This induces a $\PP^{5-s}$--bundle  $\Xcal_3(\mu,\nu)\gitquot \PGL(3)\dashrightarrow \PP\Omega_3(\mu,\nu)$. 
\end{lemma}
\begin{proof}
To see that the map is well-defined, notice that the intersection conditions of $C$ with $E$ and $Q$ give 
$$\sum_{i=1}^n\mu_ip_i+\qbf\in |\Ocal_C(3)|\textnormal{ and }\sum_{j=1}^r\nu_jq_j+\qbf\in |\Ocal_C(2)|,\textnormal{ so }\sum_{i=1}^n\mu_ip_i-\sum_{j=1}^r\nu_jq_j\sim K_C,$$
since $\omega_C\cong \Ocal_C(1)$. 
\vspace{0.5em}

Now, in order to show the map is dominant, pick a point $\left[C,(p_i)_{i=1}^n,(q_j)_{j=1}^r\right]\in \PP\Omega_3(\mu,\nu)$ such that $C$ is nonhyperelliptic and the points $q_j$ are generic on $C$. The canonical map presents the curve as a smooth plane quartic $C\in |\Ocal_{\PP^2}(4)|$. Since the points $q_j$ are generic on $C$, it follows that there exists a $\PP^{5-s}$ parametrising conics $Q$ such that $$C\cdot Q=\sum\limits_{j=1}^r\nu_jq_j+\qbf,$$
where $\qbf\in C_{8-s}$ is an effective divisor. The relation $\sum\limits_{i=1}^n\mu_ip_i-\sum\limits_{j=1}^r\nu_jq_j\sim K_C$ implies that $\sum\limits_{i=1}^n\mu_ip_i+\qbf\in |\Ocal_C(3)|$, so there exists a unique plane cubic $E$ satisfying 
$$C\cdot E=\sum_{i=1}^n\mu_ip_i+\qbf,$$
thanks to the fact that the restriction map $\H^0(\PP^2,\Ocal_{\PP^2}(3))\to \H^0(C,\Ocal_C(3))$ is an isomorphism. Finally, since $\length(E\cap Q)=6$, there exists an effective divisor $\rbf\in E_{s-2}$ such that $E\cdot Q=\qbf+\rbf$, which concludes the proof of the dominance. 
\end{proof}

At the same time, the following result establishes the rationality of $\Xcal_3(\mu,\nu)$, regardless of the type of $\mu$, as long as $\nu$ is as in (\ref{values of nu for irred strata in sum 3,4}). To achieve this, consider the variety 
$$\Ycal_3(\mu,\nu)=\left\{\left(E,Q,(p_i)_{i=1}^n,(q_j)_{j=1}^{r},\qbf\right):\begin{matrix}
    E\in |\Ocal_{\PP^2}(3)|, \textnormal{ }Q\in |\Ocal_{\PP^2}(2)|,\textnormal{ }p_i\in E,\textnormal{ }q_j\in Q,\\
    \qbf\in E_{8-s}\textnormal{, }E\cdot Q\ge \qbf\textnormal{, }\sum\limits_{i=1}^n\mu_ip_i+\qbf\in \left\lvert\Ocal_E(4)\right\rvert\\
    \textnormal{and }\sum\limits_{i=1}^n\mu_ip_i+\sum\limits_{j=1}^r\nu_jq_j+2\qbf\in |\Ocal_{E\cup Q}(4)|   
\end{matrix}\right\}$$
and the map $\pi:\Xcal_3(\mu,\nu)\dashrightarrow \Ycal_3(\mu,\nu)$ forgetting the curve $C$.

\begin{remark}\label{remark on gluing sections}
In the definition of $\Ycal_3(\mu,\nu)$, the presence of the condition $\sum\limits_{i=1}^n\mu_ip_i+\qbf\in |\Ocal_E(4)|$ is already expected, since similar requirements have appeared throughout the entire text. However, in the context at hand, the second condition on the points is also necessary, as it ensures that the section $s_E$ having divisor of zeroes $(s_E)_0=\sum\limits_{i=1}^n\mu_ip_i+\qbf$ and the section $s_Q$ having divisor of zeroes $(s_Q)_0=\sum\limits_{j=1}^r\nu_jq_j+\qbf$ also glue along the remaining intersection divisor $\rbf$ of $E$ and $Q$, in order to lift to a section $s\in \H^0(\PP^2,\Ocal_{\PP^2}(4))$. 
\vspace{0.5em}

Nevertheless, when $s=3$, the divisor $\rbf$ is a single point, so this condition is vacuous, as the agreement of $s_E$ and $s_Q$ at $\rbf$ is achieved directly by scaling. On the other hand, for $s=4$ this is no longer automatic and the condition needs to be stated. 
\end{remark}

Having said these, we are ready to state and prove the second lemma. 

\begin{lemma}
The map $\pi:\Xcal_3(\mu,\nu)\dashrightarrow \Ycal_3(\mu,\nu)$ is a birational morphism and $\Ycal_3(\mu,\nu)$ is an irreducible rational variety, whenever $\nu$ is as in (\ref{values of nu for irred strata in sum 3,4}). Thus $\Xcal_3(\mu,\nu)$ is also rational and irreducible, whenever $\nu$ satisfies (\ref{values of nu for irred strata in sum 3,4}). 
\end{lemma}
\begin{proof}
We first prove that $\pi$ is generically injective, using Theorem \ref{big indep lemma}. The divisor of conditions is $Z=\sum\limits_{i=1}^n\mu_ip_i+\sum\limits_{j=1}^r\nu_jq_j+\qbf$ and it equals $Z=3K_C+\sum\limits_{j=1}^r\nu_jq_j$ on $C$. Therefore, the number of independent conditions imposed on $|\Ocal_{\PP^2}(4)|\cong \PP^{14}$ is
\[
\length(Z)-h^0\left(C,Z-3K_C\right)=12+s-(s-2)=14,
\]
as expected. The details that remain to be checked are:
\begin{enumerate}[label=\textbf{(\roman*)}]
    \item\label{Ycal3 is rational} the variety $\Ycal_3(\mu,\nu)$ is rational, and
    \item\label{Xcal3 and Ycal3 have same dim} the varieties $\Xcal_3(\mu,\nu)$ and $\Ycal_3(\mu,\nu)$ have the same dimension, so $\pi:\Xcal_3(\mu,\nu)\dashrightarrow \Ycal_3(\mu,\nu)$ is a birational isomorphism.
\end{enumerate}
We first check \ref{Ycal3 is rational}. Introduce the variety 
$$\Zcal_3(\mu,\nu)=\left\{\left(E,(p_i)_{i=1}^n\right):E\textnormal{ plane cubic, }p_i\in E\right\},$$
together with the forgetful map $\pi_1:\Ycal_3(\mu,\nu)\dashrightarrow \Zcal_3(\mu,\nu)$. It follows immediately that $\Zcal_3(\mu,\nu)$ is rational, since $n\le 8$, and at the same time that $\pi_1$ is a dominant map, obtained as a tower of projective bundles: $\qbf$ varies in a $\PP^{7-s}$ according to the condition $\sum\limits_{i=1}^n\mu_ip_i+\qbf\in |\Ocal_E(4)|$, and $Q$ varies in a $\PP^{s-3}$--bundle according to the condition $E\cdot Q\ge \qbf$, while the points $q_j$ are retrieved as follows: 

\begin{itemize}
    \item if $s=3$, then $q_j\in Q$ can be any $r$ generic points, since the gluing conditions are vacuous, as explained in Remark \ref{remark on gluing sections}, and 
    \item if $s=4$, observe that $\nu_r=1$ for any $\nu$ as in (\ref{values of nu for irred strata in sum 3,4}); then $(q_j)_{j=1}^{r-1}$ can be any $r-1$ generic points on $Q$, and the gluing conditions recover $q_r$ uniquely.
\end{itemize}
We are left to check \ref{Xcal3 and Ycal3 have same dim}. This is a simple computation: Lemma \ref{dominance in genus 3 for odd nu 3, 4} shows that the map $\Xcal_3(\mu,\nu)/\PGL(3)\dashrightarrow \PP\Omega_3(\mu,\nu)$ is a $\PP^{5-s}$--bundle, so $\dim \Xcal_3(\mu,\nu)=16+n+r-s$. At the same time, the reconstruction of $\Ycal_3(\mu,\nu)$ from above shows 
$$\dim \Ycal_3(\mu,\nu)=\dim \Zcal_3(\mu,\nu)+(7-s)+(s-3)+(r-s+3)=16+n+r-s.$$
These show $\dim \Xcal_3(\mu,\nu)=\dim \Ycal_3(\mu,\nu)$, concluding the proof of the statement \ref{Xcal3 and Ycal3 have same dim}. Together with \ref{Ycal3 is rational}, it follows that $\Xcal_3(\mu,\nu)$ is an irreducible rational variety and the second lemma is proven. 
\end{proof}



\section{Meromorphic Differential Strata in Genus \texorpdfstring{$2$}{}}

This section handles the meromorphic differential strata in genus $2$. Theorem \ref{thm: unirat diff strata in genus 2} states that the components of 
$$\PP\Omega_2(\mu,\nu)=\left\{\left[C,(p_i)_{i=1}^n,(q_j)_{j=1}^r\right]\in \Mcal_{2,n+r}:\Ocal_C\left(\sum_{i=1}^n\mu_ip_i-\sum_{j=1}^r\nu_jq_j\right)\cong \omega_C\right\}$$
are unirational, whenever $|\nu|\le 5$ or $|\nu|=6,7$ in certain cases. The main difference between genus $2$ and genus $3$, $4$ or $5$ is that the canonical bundle is not very ample. This forces other kinds of choices, and the models that we use for a genus $2$ curve depend on the type of the partitions $\mu$ and $\nu$. Therefore, the analysis must be split into subsections, according to the value of $s=|\nu|=\sum\limits_{j=1}^r\nu_j$. 

\subsection{Meromorphic differential strata with $2$ poles}

For $\nu\in \{(2),(1,1)\}$, we prefer to write $q^{\pm}$ for the two poles, just as in the previous sections. The space of interest is
$$\PP\Omega_2(\mu,\nu)=\left\{\left[C,(p_i)_{i=1}^n,q^{\pm}\right]\in \Mcal_{2,n+r}:\Ocal_C\left(\sum_{i=1}^n\mu_ip_i-q^+-q^-\right)\cong \omega_C\right\},$$
with the same notational convention as above, that is $q^+=q^-=:q$ when $\nu=(2)$. When $\mu$ is not even, the stratum $\PP\Omega_2(\mu,\nu)$ is irreducible. When $\mu$ is even, there are two components $\PP\Omega_2(\mu,\nu)^+$ and $\PP\Omega_2(\mu,\nu)^-$. 

\vskip 0.5em
It can be easily seen that the points of the components $\PP\Omega_2(\mu,\nu)^-$ for $\nu=(2)$ and $\PP\Omega_2(\mu,\nu)^+$ for $\nu=(1,1)$ have the property that $\sum\limits_{i=1}^n\mu_i/2\cdot p_i\sim K_C$ and either $q$ is a Weierstrass point, or $(q^+,q^-)$ are involuted by the hyperelliptic structure. Therefore, these are the anti-invariant hyperelliptic components and they are known to be unirational (see e.g. \cite{chencostantinimoeller2024}, Section 2.3). For the rest of this subsection, we treat together the remaining situations, namely 
\begin{enumerate}[label=\textbf{(\alph*)}]
    \item\label{1 nontrivial case gen 2 & 2 poles} $\PP\Omega_2(\mu,\nu)^+$ for $\mu$ even, $\nu=(2)$ and $\PP\Omega_2(\mu,\nu)^-$ for $\mu$ even, $\nu=(1,1)$, as well as 
    \item\label{2 nontrivial case gen 2 & 2 poles} the irreducible strata $\PP\Omega_2(\mu,\nu)$, when $\mu$ is not even and $\nu\in \{(2),(1,1)\}$. 
\end{enumerate}
Let now $(\mu,\nu)$ be as in \ref{1 nontrivial case gen 2 & 2 poles} or \ref{2 nontrivial case gen 2 & 2 poles} and pick a point $\left[C,(p_i)_{i=1}^n,q^{\pm}\right]$ in the corresponding component. This means, in particular, that $q^++q^-\not \sim K_C$. We use this detail when we show the dominance of our desired component by the following incidence variety:

\begin{definition}
For $(\mu,\nu)$ as in \ref{1 nontrivial case gen 2 & 2 poles} or \ref{2 nontrivial case gen 2 & 2 poles}, consider the incidence variety 
$$\Xcal_2(\mu,\nu)=\left\{\left(\Gamma,L,(p_i)_{i=1}^n,q,\ell^{\pm}\right):\begin{matrix}
    \Gamma\in |\Ocal_{\PP^2}(4)|\textnormal{ nodal at }q\textnormal{ with tangent cone }\ell^{\pm},\\
    L \in (\PP^2)^\vee\textnormal{ a line with }\Gamma\cdot L=\sum\limits_{i=1}^n\mu_ip_i
\end{matrix}\right\}.$$
As in the definitions of an incidence variety for meromorphic strata of differentials in genera $3$, $4$ and $5$ with $2$ poles, the notation stands for a single line $\ell^{\pm}=:\ell$, when $\nu=(2)$. In particular, if $\nu=(2)$, then $q$ is a cusp of $\Gamma$. 
\end{definition}

Having this introduced, we proceed with the two lemmas that provide the unirationality of all components listed at \ref{1 nontrivial case gen 2 & 2 poles} and \ref{2 nontrivial case gen 2 & 2 poles}. 

\begin{lemma}\label{dominance in gen 2 for sum(nu)=2}
The maps $\Xcal_2(\mu,\nu)\dashrightarrow \PP\Omega_2(\mu,\nu)^+$ for $\mu$ even and $\nu=(2)$, $\Xcal_2(\mu,\nu)\dashrightarrow \PP\Omega_2(\mu,\nu)^-$ for $\mu$ even and $\nu=(1,1)$, and $\Xcal_2(\mu,\nu)\dashrightarrow\PP\Omega_2(\mu,\nu)$ for $\mu$ not even, that are given by 
$$\left[\Gamma,L,(p_i)_{i=1}^n,q,\ell^{\pm}\right]\mapsto \left[C,(p_i)_{i=1}^n,q^{\pm}\right],$$
are well-defined and dominant, where $\sigma:C\to \Gamma$ is the normalisation map, and $q^{\pm}\in C$ correspond to the tangent directions $\ell^{\pm}$ of $\Gamma$ at $q$. Quotienting by $\PGL(3)$, we obtain a birational isomorphism $\Xcal_3(\mu, \nu)\gitquot\PGL(3) \dashrightarrow\PP\Omega_2(\mu, \nu)$ (respectively $\Xcal_3(\mu, \nu)\gitquot\PGL(3) \dashrightarrow\PP\Omega_2(\mu, \nu)^\pm$).
\end{lemma}
\begin{proof}
To see that these maps are well-defined, the condition $\Gamma\cdot L=\sum\limits_{i=1}^n\mu_ip_i$ gives 
$$\omega_C(q^++q^-)=\sigma^*\omega_\Gamma=\sigma^*\Ocal_\Gamma(1)=\Ocal_C\left(\sum_{i=1}^n\mu_ip_i\right),$$
whence $\Ocal_C\left(\sum\limits_{i=1}^n\mu_ip_i-q^+-q^-\right)\cong \omega_C$. As for the dominance statement, start with a point $\left[C,(p_i)_{i=1}^n,q^{\pm}\right]$ in the corresponding component of $\PP\Omega_2(\mu,\nu)$, so $\sum\limits_{i=1}^n\mu_ip_i\sim K_C+q^++q^-$ and $q^++q^-\not\sim K_C$. The line bundle $A:=\omega_C(q^++q^-)$ is a $\gfrak^2_4$ and maps $C$ birationally to a plane quartic $\Gamma\in |\Ocal_{\PP^2}(4)|$, having a node below $(q^+,q^-)$ if $\nu=(1,1)$, or a cusp at $q$, if $\nu=(2)$. Finally, since $\sum\limits_{i=1}^n\mu_ip_i\in |A|$, there exists a line $L\subset \PP^2$ such that $C\cdot L=\sum\limits_{i=1}^n\mu_ip_i$. The dominance claim is now established. 
\end{proof}

We close this subsection with the following lemma. 

\begin{lemma}
The variety $\Xcal_2(\mu,\nu)$ is rational, whenever $\mu$ and $\nu$ are as in \ref{1 nontrivial case gen 2 & 2 poles} or \ref{2 nontrivial case gen 2 & 2 poles}. 
\end{lemma}
\begin{proof}
Consider the variety 
$$\Ycal_2(\mu,\nu)=\left\{\left(L,(p_i)_{i=1}^n\right):L\textnormal{ line, }p_i\in L\right\}\times \left\{\left(q,\ell^{\pm}\right):q\in \PP^2\textnormal{, }\ell^{\pm}\textnormal{ lines through }q\right\}$$
which is clearly rational. There is a forgetful map $\pi:\Xcal_2(\mu,\nu)\dashrightarrow \Ycal_2(\mu,\nu)$ that remembers all data, but the curve $\Gamma$. The following two claims show that $\Xcal_2(\mu,\nu)$ is indeed rational:

\begin{enumerate}[label=\textbf{(\roman*)}]
    \item\label{indep lemma gen 2 & 2 poles} the map $\pi$ is a $\PP^5$--bundle over the image, and 
    \item\label{dim count gen 2 & 2 poles} $\dim \Xcal_2(\mu,\nu)=9+n+r$ and $\dim \Ycal_2(\mu,\nu)=4+n+r$, so $\pi$ is, in fact, a $\PP^5$--bundle over a dense open in the target. 
\end{enumerate}

For \ref{indep lemma gen 2 & 2 poles}, we use Theorem \ref{big indep lemma}. We look on the blow-up $S:=\textnormal{Bl}_q\left(\PP^2\right)$ of $\PP^2$ at $q$ and we denote the normalisation of $\Gamma$ by $C\in |4H-2E_q|$. Then the divisor of conditions is $Z=\sum\limits_{i=1}^n\mu_ip_i+q^++q^-$ and on $C$ it reads $Z=K_C+2(q^++q^-)$, so the number of independent conditions imposed on $|4H-2E_q|\cong \PP^{11}$ is $\length(Z)-h^0(C,Z-K_C-2H\vert_C)=6$. Regarding \ref{dim count gen 2 & 2 poles}, simply notice that Lemma \ref{dominance in gen 2 for sum(nu)=2} gives a birational isomorphism from $\Xcal_2(\mu,\nu)\gitquot\PGL(3)$ to the corresponding component of $\PP\Omega_2(\mu,\nu)$. The second lemma is now proven. 
\end{proof}

\subsection{Meromorphic differential strata in genus $2$ with $3$ poles}

When the type of poles is $\nu\in\{(1,1,1),(1,2),(3)\}$, the meromorphic strata of differentials $\PP\Omega_2(\mu,\nu)$ is irreducible. The unirational parametrisation of $\PP\Omega_2(\mu,\nu)$ is given by 

\begin{definition}
Consider the surface $S=\PP^1\times \PP^1$. For $|\nu|=3$, consider the incidence variety 
$$\Xcal_2(\mu,\nu)=\left\{\left(C,Q,\ell,(p_i)_{i=1}^n,(q_j)_{j=1}^r\right):\begin{matrix}
    C\in |\Ocal_S(3,2)|\textnormal{, }Q\in |\Ocal_S(1,1)|\textnormal{, }\ell\in |\Ocal_S(0,1)|\\
    \textnormal{with }C\cdot Q=\sum\limits_{i=1}^n\mu_ip_i\textnormal{ and }C\cdot \ell=\sum\limits_{j=1}^r\nu_jq_j
\end{matrix}\right\}.$$
\end{definition}

\begin{lemma}\label{dominance in gen 2 for sum(nu)=3}
The map $\Xcal_2(\mu,\nu)\dashrightarrow \PP\Omega_2(\mu,\nu)$ given by 
$$\left[C,Q,\ell,(p_i)_{i=1}^n,(q_j)_{j=1}^r\right]\mapsto \left[C,(p_i)_{i=1}^n,(q_j)_{j=1}^r\right]$$
is well-defined and dominant. Moreover, the induced map $\Xcal_2(\mu,\nu) \gitquot \PGL(2)\times \PGL(2) \dashrightarrow \PP\Omega_2(\mu, \nu)$ is a birational isomorphism.
\end{lemma}
\begin{proof}
To check that the map is indeed well-defined, notice that the adjunction formula gives $\omega_C=\Ocal_C(1,0)$ for a curve $C\in |\Ocal_S(3,2)|$, so if $\left(C,Q,\ell,(p_i)_{i=1}^n,(q_j)_{j=1}^r,\right)$ is a point of $\Xcal_2(\mu,\nu)$, then the intersection conditions of $C$ with $Q$ and $\ell$ imply 
$$\sum_{i=1}^n\mu_ip_i\in |\Ocal_C(1,1)|\textnormal{ and }\sum_{j=1}^r\nu_jq_j\in |\Ocal_C(0,1)|,\textnormal{ so }\sum_{i=1}^n\mu_ip_i-\sum_{j=1}^r\nu_jq_j\sim K_C,$$
as expected. As for the dominance claim, start with a point $\left[C,(p_i)_{i=1}^n,(q_j)_{j=1}^r\right]\in \PP\Omega_2(\mu,\nu)$ such that the points $q_j$ are in general position on $C$. Then the line bundle $A:=\Ocal_C\left(\sum\limits_{j=1}^r\nu_jq_j\right)$ is a basepoint free $\gfrak^1_3$ on $C$ and, together with the hyperelliptic structure $\gfrak^1_2$ on $C$, gives an embedding $C\hookrightarrow S$, such that $C\in |\Ocal_S(3,2)|$. Furthermore, there exists $\ell\in |\Ocal_S(0,1)|$ such that $C\cdot \ell=\sum\limits_{j=1}^r\nu_jq_j$, and then the condition $\sum\limits_{i=1}^n\mu_ip_i-\sum\limits_{j=1}^r\nu_jq_j\sim K_C$ implies
$$\sum_{i=1}^n\mu_ip_i\in |\Ocal_C(1,1)|,\textnormal{ so there exists a unique }Q\in |\Ocal_S(1,1)|\textnormal{ with }C\cdot Q=\sum_{i=1}^n\mu_ip_i.$$
This concludes the first lemma. 
\end{proof}

Now we continue with the second lemma, which covers all possible partitions $\mu$ of $5$. 

\begin{lemma}
The variety $\Xcal_2(\mu,\nu)$ is rational, whenever $\nu$ has $|\nu|=3$. 
\end{lemma}
\begin{proof}
Consider the variety 
$$\Ycal_2(\mu,\nu)=\left\{\left(Q,(p_i)_{i=1}^n\right):Q\in |\Ocal_S(1,1)|\textnormal{, }p_i\in Q\right\}\times \left\{\left(\ell,(q_j)_{j=1}^r\right):\ell\in |\Ocal_S(0,1)|\textnormal{, }q_j\in \ell\right\}$$
which is clearly rational. There is a forgetful map $\pi:\Xcal_2(\mu,\nu)\dashrightarrow \Ycal_2(\mu,\nu)$ that remembers all data, except the curve $C$. We check the following two claims, which are enough to conclude that $\Xcal_2(\mu,\nu)$ is rational: 
\begin{enumerate}[label=\textbf{(\roman*)}]
    \item\label{indep lemma gen 2 & 3 poles} the map $\pi$ is a $\PP^3$--bundle over the image, and 
    \item\label{dim count gen 2 & 3 poles} $\dim\Xcal_2(\mu,\nu)=7+n+r$ and $\dim \Ycal_2(\mu,\nu)=4+n+r$, so $\pi$ is a $\PP^3$--bundle over a dense open in the target. 
\end{enumerate}

The \ref{indep lemma gen 2 & 3 poles} part follows from Theorem \ref{big indep lemma}. Indeed, the divisor of conditions is $Z=\sum\limits_{i=1}^n\mu_ip_i+\sum\limits_{j=1}^r\nu_jq_j$, so the number of independent conditions imposed on $|\Ocal_S(3,2)|\cong \PP^{11}$ is $\length(Z)-h^0(C,Z-2K_C-2\Ocal_C(0,1))=8$. In terms of \ref{dim count gen 2 & 3 poles}, by Lemma \ref{dominance in gen 2 for sum(nu)=3} it follows that $\dim \Xcal_2(\mu, \nu) = 7 + n + r$. The second lemma is therefore proven and this concludes the situation of $|\nu|=3$ poles. 
\end{proof}

\subsection{Meromorphic differential strata in genus $2$ with $4$ poles}

In this case, we need to split the analysis since the space $\PP\Omega_2(\mu,\nu)$ is reducible, for certain tuples $(\mu,\nu)$. Recall that, when $\mu$ and $\nu$ are even, the components of $\PP\Omega_2(\mu, \nu)$ are $\PP\Omega_2(\mu,\nu)^+$ and $\PP\Omega_2(\mu,\nu)^-$.
\vspace{0.5em}

We develop the analysis of meromorphic differential strata with $4$ poles into the following parts:

\begin{enumerate}[label=\textbf{(\alph*)}]
    \item\label{case mu nu even} \textit{Case of $\mu$ and $\nu$ even,} when there are two components $\PP\Omega_2(\mu,\nu)^+$ and $\PP\Omega_2(\mu,\nu)^-$, 
    \item\label{case other mu nu, for sum(nu)=4} \textit{Case of $\mu_n=1$ or $\nu_r=1$,} corresponding to one irreducible component $\PP\Omega_2(\mu,\nu)$, 
    \item\label{case other mu nu, for sum(nu)=4, bis} \textit{Case of $\mu=(3,3)$ and $\nu$ even}, which still corresponds to one irreducible component $\PP\Omega_2(\mu,\nu)$. 
\end{enumerate}

For the moment, we treat \ref{case mu nu even}, \ref{case other mu nu, for sum(nu)=4} and \ref{case other mu nu, for sum(nu)=4, bis} simultaneously. The model for a genus $2$ curve in these cases is derived from the $4$ poles: if the $q_j$ are general, then $A:=\Ocal_C\left(\sum\limits_{i=1}^r\nu_jq_j\right)$ is a $\gfrak^2_4$ that maps the curve to a singular plane quartic in $\PP^2$. This suggests the following 

\begin{definition}\label{incidence var other mu nu for sum(nu)=4}
For $\mu$ and $\nu$ as in \ref{case mu nu even}, \ref{case other mu nu, for sum(nu)=4} or \ref{case other mu nu, for sum(nu)=4, bis}, consider the incidence variety 
$$\Xcal_2(\mu,\nu)=\left\{\left(\Gamma,Q,\ell,x,(p_i)_{i=1}^n,(q_j)_{j=1}^r\right):\begin{matrix}
    \Gamma\in |\Ocal_{\PP^2}(4)|\textnormal{ nodal at }x\textnormal{, }Q\textnormal{ conic, }\ell\textnormal{ line}\\
    \textnormal{with }\Gamma\cdot Q=2x+\sum\limits_{i=1}^n\mu_ip_i\textnormal{, }\Gamma\cdot \ell=\sum\limits_{j=1}^r\nu_jq_j
\end{matrix}\right\}.$$
\end{definition}

The justification for the choice of the variety $\Xcal_2(\mu,\nu)$ lies in the following lemma:

\begin{lemma}\label{dominance in gen 2 for sum(nu)=4 other mu nu}
The map $\Xcal_2(\mu,\nu)\dashrightarrow \PP\Omega_2(\mu,\nu)$ given by 
$$\left[\Gamma,Q,\ell,x,(p_i)_{i=1}^n,(q_j)_{j=1}^r\right]\mapsto \left[C,(p_i)_{i=1}^n,(q_j)_{j=1}^r\right]$$
is well-defined and dominant over each component of the target, where $\sigma:C\to\Gamma$ is the normalisation map. Moreover, the map $\Xcal_2(\mu, \nu) \gitquot \PGL(3) \dashrightarrow \PP\Omega_2(\mu, \nu)$ is a birational isomorphism.
\end{lemma}
\begin{proof}
First, we show that the map is well-defined. For a point $\left(\Gamma,Q,\ell,x,(p_i)_{i=1}^n,(q_j)_{j=1}^r\right)$ of $\Xcal_2(\mu,\nu)$ and the normalisation $\sigma:C\to \Gamma$, let $\sigma^{-1}(x)=\{x^{\pm}\}$ and we have 
$$\omega_C\left(x^++x^-\right)=\sigma^*\omega_\Gamma=\sigma^*\Ocal_{\Gamma}(1)=\sigma^*\Ocal_{\Gamma}(Q-\ell)=\Ocal_C\left(x^++x^-+\sum_{i=1}^n\mu_ip_i-\sum_{j=1}^r\nu_jq_j\right),$$
whence $\sum\limits_{i=1}^n\mu_ip_i-\sum\limits_{j=1}^r\nu_jq_j\sim K_C$. For the dominance claim, pick $\left[C,(p_i)_{i=1}^n,(q_j)_{j=1}^r\right]$ in $\PP\Omega_2(\mu,\nu)$ such that the $q_j$ are generic. As mentioned in the paragraph just above Definition \ref{incidence var other mu nu for sum(nu)=4}, the line bundle $A=\Ocal_C\left(\sum\limits_{j=1}^r\nu_jq_j\right)$ sends the curve $C$ to a plane quartic $\Gamma\in|\Ocal_{\PP^2}(4)|$ which is singular at some point $x\in \PP^2$. There also exists a line $\ell$ that cuts precisely the divisor $\sum\limits_{j=1}^r\nu_jq_j$ on $C$. With these in mind, let $\{x^{\pm}\}$ be the preimage of $x$, so we can write 
$$\sum_{i=1}^n\mu_ip_i+x^++x^-\sim K_C\left(\sum_{j=1}^r\nu_jq_j+x^++x^-\right)\sim \sigma^*\Ocal_{\Gamma}(2),$$
and thus there exists a unique conic $Q$ such that $\Gamma\cdot Q=2x+\sum\limits_{i=1}^n\mu_ip_i$. This establishes the first lemma. 
\end{proof}

As expected, there is now a second lemma, that treats the rationality of $\Xcal_2(\mu,\nu)$. 

\begin{lemma}
Let $\mu$ and $\nu$ be as in \ref{case mu nu even}, \ref{case other mu nu, for sum(nu)=4} or \ref{case other mu nu, for sum(nu)=4, bis}. All components of the variety $\Xcal_2(\mu,\nu)$ are rational. In particular, all components of $\PP\Omega_2(\mu,\nu)$ are unirational.
\end{lemma}
\begin{proof}
Define the variety 
$$\Ycal_2(\mu,\nu)=\left\{\left(Q,\ell,x,(p_i)_{i=1}^n,(q_j)_{j=1}^r\right):
\begin{matrix}
    Q\textnormal{ conic, }\ell\textnormal{ line, }x,p_i\in Q\textnormal{, }q_j\in \ell,\\
    2x+\sum\limits_{i=1}^n\mu_ip_i+\sum\limits_{j=1}^r\nu_jq_j\in |\Ocal_{Q\cup \ell}(4)|
\end{matrix}\right\}$$
and consider the forgetful map $\pi:\Xcal_2(\mu,\nu)\dashrightarrow\Ycal_2(\mu,\nu)$ that remembers all the data, except the curve $\Gamma$. The following three claims are sufficient, in order to conclude that all components of $\Xcal_2(\mu,\nu)$ are rational varieties:
\begin{enumerate}[label=\textbf{(\roman*)}]
    \item\label{indep lemma gen 2 & 4 poles other mu nu} the map $\pi$ is a $\PP^2$--bundle over the image, 
    \item\label{rat check for Ycal2} $\Ycal_2(\mu,\nu)$ has at most two components and they are rational, whenever $\mu$, $\nu$ are as in \ref{case mu nu even}, \ref{case other mu nu, for sum(nu)=4} or \ref{case other mu nu, for sum(nu)=4, bis}, and
    \item\label{dim count gen 2 & 4 poles other mu nu} $\dim \Xcal_2(\mu,\nu)=9+n+r$ and $\dim\Ycal_2(\mu,\nu)=7+n+r$, so $\pi$ is, in fact, a $\PP^2$--bundle over dense opens in each component of the target. 
\end{enumerate}

We use Theorem \ref{big indep lemma} to check \ref{indep lemma gen 2 & 4 poles other mu nu}. Let $S=\textnormal{Bl}_x\left(\PP^2\right)$ be the blow-up of $\PP^2$ at $x$ and let $C\in |4H-2E_x|$ be the strict transform of $\Gamma$. The divisor of conditions is $Z=\sum\limits_{i=1}^n\mu_ip_i+\sum\limits_{j=1}^r\nu_jq_j$, which on $C$ is equal to $Z=K_C+2\sum\limits_{j=1}^r\nu_jq_j$, so the number of independent conditions imposed on $|4H-2E_x|\cong \PP^{11}$ is $\length(Z)-h^0\left(C,Z-K_C-2H\vert_C\right)=10-h^0(C,\Ocal_C)=9$. 
\vspace{0.5em}

Next, let us address part \ref{rat check for Ycal2}, by separating the argument according to the cases \ref{case mu nu even}, \ref{case other mu nu, for sum(nu)=4} and \ref{case other mu nu, for sum(nu)=4, bis}. 
\vspace{0.5em}

Take first the situation at \ref{case mu nu even}, $\mu$ and $\nu$ are both even. The condition on points becomes $x+\sum\limits_{i=1}^n\mu_i/2\cdot p_i+\sum\limits_{j=1}^r\nu_j/2\cdot q_j\in |\Ocal_{Q\cup \ell}(2)\otimes \tau|$, for one of the two $2$-torsion line bundle $\tau$ on $Q\cup \ell$.

If $\tau\cong \Ocal_{Q\cup \ell}$, then $x+\sum\limits_{i=1}^n\mu_i/2\cdot p_i+\sum\limits_{j=1}^r\nu_j/2\cdot q_j\in |\Ocal_{Q\cup \ell}(2)|$, so there exists a unique conic $Q'\ne Q$ such that $Q\cdot Q'=x+\sum\limits_{i=1}^n\mu_i/2\cdot p_i$ and $\ell\cdot Q'=\sum\limits_{j=1}^r\nu_j/2\cdot q_j$. These imply $\Gamma\cdot Q'\ge 2x+\sum\limits_{i=1}^n\mu_i/2\cdot p_i$ and, by comparing with the relation $\Gamma\cdot Q=2x+\sum\limits_{i=1}^n\mu_ip_i$, we deduce that $\sum\limits_{i=1}^n\mu_i/2\cdot p_i$ is a moving divisor on $\Gamma$. Therefore, the situation $\tau\cong \Ocal_{Q\cup \ell}$ corresponds to the component $\PP\Omega_2(\mu,\nu)^-$, whereas the situation $\tau\not \cong \Ocal_{Q\cup \ell}$ corresponds to the other component $\PP\Omega_2(\mu,\nu)^+$. The components of $\Ycal_2(\mu,\nu)$ when $\mu$ and $\nu$ are even are thus described by 
$$\Ycal_2(\mu,\nu)^-:=\left\{\left(Q,\ell,x,(p_i)_{i=1}^n,(q_j)_{j=1}^r\right):\begin{matrix}
    Q\textnormal{ conic, }\ell\textnormal{ line, }x,p_i\in Q,\textnormal{ }q_j\in \ell\\
    x+\sum\limits_{i=1}^n\mu_i/2\cdot p_i+\sum\limits_{j=1}^r\nu_j/2\cdot q_j\in |\Ocal_{Q\cup \ell}(2)|
\end{matrix}\right\}$$
and 
$$\Ycal_2(\mu,\nu)^+:=\left\{\left(Q,\ell,x,(p_i)_{i=1}^n,(q_j)_{j=1}^r\right):\begin{matrix}
    Q\textnormal{ conic, }\ell\textnormal{ line, }x,p_i\in Q\textnormal{, }q_j\in\ell,\\
    x+\sum\limits_{i=1}^n\mu_i/2\cdot p_i+\sum\limits_{j=1}^r\nu_j/2\cdot q_j\in |\Ocal_{Q\cup \ell}(2)\otimes \tau|\\
    \textnormal{ for the unique }\tau\ne \Ocal_{Q\cup \ell}\textnormal{ with }\tau^2\cong \Ocal_{Q\cup \ell}
\end{matrix}\right\}$$
There are $\PP^2$--bundles 
$$\pi^+:\Xcal_2(\mu,\nu)^+\dashrightarrow \Ycal_2(\mu,\nu)^+\textnormal{ and }\pi^-:\Xcal_2(\mu,\nu)^-\dashrightarrow \Ycal_2(\mu,\nu)^-,$$ 
where $\Xcal_2(\mu,\nu)^+$ and $\Xcal_2(\mu,\nu)^-$ are the corresponding components of $\Xcal_2(\mu,\nu)$. These show that, in any case, the point $x$ is uniquely determined in terms of the $p_i$ and the $q_j$. Therefore, both components of $\Ycal_2(\mu,\nu)$ are birationally isomorphic to 
$$\left\{(Q,(p_i)_{i=1}^n):Q\textnormal{ conic, }p_i\in Q\right\}\times \left\{(\ell,(q_j)_{j=1}^r):\ell\textnormal{ line, }q_j\in \ell\right\},$$
and the latter variety is clearly rational. For the case \ref{case other mu nu, for sum(nu)=4}, if $\mu_n=1$, then $\Ycal_2(\mu,\nu)$ is irreducible, and $p_n$ is uniquely determined by the points $x,(p_i)_{i=1}^{n-1}$ and $q_j$, so in this case $\Ycal_2(\mu,\nu)$ is birationally isomorphic to 
$$\left\{(Q,x,(p_i)_{i=1}^{n-1}):Q\textnormal{ conic, }x,p_i\in Q\right\}\times \left\{(\ell,(q_j)_{j=1}^r):\ell\textnormal{ line, }q_j\in \ell\right\},$$
clearly a rational variety. The subcase $\nu_r=1$ is the same as $\mu_n=1$, and therefore we omit it. The last case \ref{case other mu nu, for sum(nu)=4, bis} is concerned with the partitions $\mu=(3,3)$ and $\nu\in \{(2,2),(4)\}$. If $\nu=(2,2)$ and $\left(Q,\ell,x,(p_i)_{i=1}^2,{q_j}_{j=1}^2\right)\in \Ycal_2(\mu,\nu)$, we write $p_1+p_2+q_1+q_2\in |\Ocal_{Q\cup \ell}(1)(y)|$, for a unique point $y\in \ell$. The condition $x+3p_1+3p_2+2q_1+2q_2\in |\Ocal_{Q\cup \ell}(4)|$ on points then becomes $p_1+p_2+2x+2y\in |\Ocal_{Q\cup \ell}(2)|$, and it is now clear that $\Ycal_2(\mu,\nu)$ is birationally isomorphic to the variety 
$$\left\{\left(Q,\ell,x,y,p_2,q_2\right):Q\textnormal{ conic, }\ell\textnormal{ line, }p_2,x\in Q\textnormal{, }q_2,y\in \ell\right\},$$
which is, of course, rational. One deals in a similar manner with the subcase $\mu=(3,3)$, $\nu=(4)$. Finally, the \ref{dim count gen 2 & 4 poles other mu nu} part follows from Lemma \ref{dominance in gen 2 for sum(nu)=4 other mu nu}, which shows that there is a birational isomorphism $\Xcal_2(\mu,\nu)\gitquot\PGL(3)\dashrightarrow\PP\Omega_2(\mu,\nu)$, whilst the dimension estimate for $\Ycal_2(\mu,\nu)$ is an automatic consequence of the discussion around \ref{rat check for Ycal2}. 
\end{proof}

\subsection{Meromorphic strata with $5$, $6$ or $7$ poles}

When the type $\nu$ of poles satisfies $|\nu|=5$ or $|\nu|=7$, the meromorphic strata $\PP\Omega_2(\mu,\nu)$ are irreducible. When $|\nu|=6$, the stratum $\PP\Omega_2(\mu,\nu)$ has two components $\PP\Omega_2(\mu,\nu)^+$ and $\PP\Omega_2(\mu,\nu)^-$.
\vspace{0.5em}

Set $s=|\nu|\in \{5,6,7\}$. For the particular class of partitions of poles described below at \ref{case of at most 4 distinct poles}, we offer a unified argument, which is independent of whether $s$ is $5$, $6$ or $7$. 

\begin{enumerate}[label=\textbf{(\alph*)}]
    \item \label{case of at most 4 distinct poles} \textit{Case of at most $4$ distinct poles,} where we further require that $\nu$ can be written as $\nu=2\nu'+\nu''$, such that $\nu'$, $\nu''$ are tuples of nonnegative integers having $|\nu'|=s-4$, except $(6,1)$, $(4,3)$, $(7)$.
    \item\label{case other mu nu, for sum(nu)=5,6} \textit{Case of $\nu=(1,1,1,1,1)$ and any other $\mu$ with $|\mu|=7$}.  
\end{enumerate}
We recall that Theorem \ref{thm: unirat diff strata in genus 2} is concerned only with the case $n\le 7$ when $s=6$ and $\nu$ is not even, and only with $n\le 8$ for the partitions $\nu$ with $s=7$. Concretely, the partitions $\nu$ that are singled out in \ref{case of at most 4 distinct poles} are
\begin{align}\label{nu with sum 5 in genus 2}
    (2,1,1,1),(2,2,1),(3,1,1),(4,1),(3,2),(5)\textnormal{ for }s=5,
\end{align}
\begin{align}\label{nu with sum 6 in genus 2}
    (2,2,1,1),(4,1,1),(3,2,1),(2,2,2),(5,1),(4,2),(3,3),(6)\textnormal{ for }s=6\textnormal{, and}
\end{align}
\begin{align}\label{nu with sum 7 in genus 2}
    (2,2,2,1),(4,2,1),(3,2,2),(5,2)\textnormal{ for }s=7.
\end{align}
Notice that, for $s=5$, these leave out only the case $\nu=(1,1,1,1,1)$, which will be treated separately, as announced in \ref{case other mu nu, for sum(nu)=5,6}. For $s=6$, the partitions listed at (\ref{nu with sum 6 in genus 2}) are exactly those announced in Theorem \ref{thm: unirat diff strata in genus 2}. The model that we use for \ref{case of at most 4 distinct poles} is described in the following

\begin{definition}
Let $\mu$, $\nu$ be as in \ref{case of at most 4 distinct poles} and write $\nu=2\nu'+\nu''$, where $|\nu''|=8-s$. Consider the incidence variety 
$$\Xcal_{2}(\mu,\nu)=\left\{\left(\Gamma,E,\ell,(p_i)_{i=1}^n,(q_j)_{j=1}^r,x\right):\begin{matrix}
    \Gamma\in |\Ocal_{\PP^2}(4)|\textnormal{ nodal at }x\textnormal{, } E\in |\Ocal_{\PP^2}(3)|\\
    \textnormal{with }\Gamma\cdot E=\sum\limits_{i=1}^n\mu_ip_i+\sum\limits_{j=1}^r\nu_j''q_j+2x,\\
    \ell\in\left(\PP^2\right)^\vee\textnormal{ with }\Gamma\cdot \ell=\sum\limits_{j=1}^r(\nu_j'+\nu_j'')q_j
\end{matrix}\right\}.$$
\end{definition}

The construction of the incidence variety, as well as the idea behind considering the specific partitions at \ref{case of at most 4 distinct poles} are presented in the first lemma below. 

\begin{lemma}\label{dominance in gen 2 for at most 4 distinct poles}
The map $\Xcal_2(\mu,\nu)\dashrightarrow \PP\Omega_2(\mu,\nu)$ given by 
$$\left[\Gamma,E,\ell,(p_i)_{i=1}^n,(q_j)_{j=1}^r,x\right]\mapsto \left[C,(p_i)_{i=1}^n,(q_j)_{j=1}^r\right]$$
is well-defined and dominant over each irreducible component of the target, where we write $\sigma:C\to \Gamma$ for the normalisation map. Furhtermore, the map $\Xcal_2(\mu, \nu)\gitquot\PGL(3) \dashrightarrow \PP\Omega_2(\mu, \nu)$ is a birational isomorphism.
\end{lemma}
\begin{proof}
For a point $\left[\Gamma,E,\ell,(p_i)_{i=1}^n,(q_j)_{j=1}^r,x\right]$ of $\Xcal_2(\mu,\nu)$, let $\sigma^{-1}(x)=\left\{x^{\pm}\right\}$ and write 
$$\omega_C(x^{+}+x^{-})=\sigma^*\omega_\Gamma=\sigma^*\Ocal_\Gamma(1)=\sigma^*\Ocal_\Gamma(E-2\ell)=\Ocal_C\left(x^++x^-+\sum_{i=1}^n\mu_ip_i-\sum_{j=1}^n\nu_jq_j\right),$$
giving precisely $\sum\limits_{i=1}^n\mu_ip_i-\sum\limits_{j=1}^r\nu_jq_j\sim K_C$. This shows that the map is well-defined. For the dominance statement, pick a point $\left[C,(p_i)_{i=1}^n,(q_j)_{j=1}^r\right]$ in any component of $\PP\Omega_2(\mu,\nu)$ such that the $q_j$ are generic on $C$. The line bundle $A:=\Ocal_C\left(\sum\limits_{j=1}^r\left(\nu_j'+\nu_j''\right)q_j\right)\in W^2_4(C)$ maps $C$ to a plane quartic $\Gamma\in |\Ocal_{\PP^2}(4)|$ that is nodal at some point $x\in \PP^2$. There exists a line $\ell\in \left(\PP^2\right)^\vee$ such that $\Gamma\cdot \ell=\sum\limits_{j=1}^r(\nu_j'+\nu_j'')q_j$, and the condition $\sum\limits_{i=1}^n\mu_ip_i-\sum\limits_{j=1}^r\nu_jq_j\sim K_C$ gives 
$$\sum_{i=1}^n\mu_ip_i+\sum_{j=1}^n\nu_j''q_j+x^++x^-\sim K_C\left(2\sum\limits_{j=1}^r(\nu_j'+\nu_j'')q_j+x^++x^-\right)\cong \sigma^*\Ocal_{\Gamma}(3),$$
so there exists a unique cubic $E$ such that $\Gamma\cdot E=\sum\limits_{i=1}^n\mu_ip_i+\sum\limits_{j=1}^r\nu_j''q_j+2x$. This settles the dominance claim and the lemma is proven. 
\end{proof}

So far, we heavily used the requirement that $\nu$ can be written as $\nu=2\nu'+\nu''$ for some tuples $\nu',\nu''$ of nonnegative integers, with $|\nu''|=8-s$. The next lemma shows that, except when $\nu\in \{(6,1),(4,3),(7)\}$, these suffice in order to conclude part \ref{case of at most 4 distinct poles}. 

\begin{lemma}
For $\mu$ and $\nu$ as in \ref{case of at most 4 distinct poles}, all components of the variety $\Xcal_2(\mu,\nu)$ are rational. In particular, all components of $\PP\Omega_2(\mu,\nu)$ are unirational. 
\end{lemma}
\begin{proof}
Consider the variety 
$$\Ycal_2(\mu,\nu)=\left\{\left(E,\ell,(p_i)_{i=1}^n,(q_j)_{j=1}^r,x\right):\begin{matrix}
    E\textnormal{ cubic, }\ell\textnormal{ line with }E\cdot \ell\ge \sum\limits_{j=1}^r\nu_j''q_j\textnormal{, }q_j\in \ell,\\
    \sum\limits_{i=1}^n\mu_ip_i+\sum\limits_{j=1}^r\nu_j''q_j+2x\in |\Ocal_E(4)|\textnormal{ and}\\    \sum\limits_{i=1}^n\mu_ip_i+\sum\limits_{j=1}^r(\nu_j'+2\nu_j'')q_j+2x\in |\Ocal_{E\cup \ell}(4)|
\end{matrix}\right\}$$
together with the map $\pi:\Xcal_2(\mu,\nu)\dashrightarrow\Ycal_2(\mu,\nu)$ that forgets the curve $\Gamma$. In order to obtain that all components of $\Xcal_2(\mu,\nu)$ are rational, it suffices to check the following three claims:
\begin{enumerate}[label=\textbf{(\roman*)}]
    \item\label{indep lemma gen 2 & at most 4 distinct poles} The map $\pi$ is a $\PP^{1}$--bundle over the image for $s=5$, and generically injective for $s=6,7$,  
    \item\label{rat check for Ycal2 when at most 4 distinct poles} $\Ycal_2(\mu,\nu)$ has at most two components and they are rational, whenever $\mu$, $\nu$ are as in \ref{case of at most 4 distinct poles}, and 
    \item\label{dim count gen 2 & at most 4 distinct poles} $\dim \Xcal_2(\mu,\nu)=9+n+r$ and either $\dim \Ycal_2(\mu,\nu)=8+n+r$ for $s=5$, or $\dim \Ycal_2(\mu,\nu)=9+n+r$ for $s=6,7$, so $\pi$ is dominant over each irreducible component of $\Ycal_2(\mu,\nu)$. 
\end{enumerate}

The claim at \ref{indep lemma gen 2 & at most 4 distinct poles} follows from Theorem \ref{big indep lemma}. As argued before, we let $S=\textnormal{Bl}_x\left(\PP^2\right)$ be the blow-up of the plane at $x$, we let $C\in |4H-2E_x|$ be the strict transform of $\Gamma$ and we keep the same notation for $E$ and $\ell$ on $S$. The divisor of conditions is $Z=\sum\limits_{i=1}^n\mu_ip_i+\sum\limits_{j=1}^r(\nu_j'+\nu_j'')q_j$, so there are $$\length(Z)-h^0\left(C,Z-K_C-2H\vert_C\right)=6+s-h^0\left(\sum\limits_{j=1}^r\nu_j'q_j\right)=\begin{cases}
    10\textnormal{ for }s=5,\textnormal{ and}\\
    11\textnormal{ for }s=6,7
\end{cases}$$
independent conditions imposed on the system $|4H-2E_x|\cong \PP^{11}$. 
\vspace{0.5em}

The dimension estimate on $\Xcal_2(\mu,\nu)$ at \ref{dim count gen 2 & at most 4 distinct poles} is an immediate consequence of Lemma \ref{dominance in gen 2 for at most 4 distinct poles}, as it gives a birational isomorphism $\Xcal_2(\mu,\nu)/\PGL(3)\dashrightarrow \PP\Omega_2(\mu,\nu)$, while the dimension estimate for $\Ycal_2(\mu,\nu)$ can be tracked through the arguments below. 

\vskip 0.5em
Let us focus on the affirmation \ref{rat check for Ycal2 when at most 4 distinct poles}. This needs to be verified step by step. We begin with the quick remark that the second condition on the points in the definition of $\Ycal_2(\mu,\nu)$ is only relevant when $s=7$, since $\length(E\cap \ell)=3$ and $\deg\left(\sum\limits_{j=1}^r\nu_j''q_j\right)=8-s=1$, so there are $2$ extra points of intersection where the relevant sections on $E$ and $\ell$ need to agree, before they lift to the same degree $4$ form $s\in \H^0(\PP^2,\Ocal_{\PP^2}(4))$, as it was explained in a similar scenario in Remark \ref{remark on gluing sections}. 
\vspace{0.5em}

\textbf{\underline{Case $s=5$}.} Here, we actually have $E\cdot \ell=\sum\limits_{j=1}^r\nu_j''q_j$, and thus the only relevant condition on points is equivalent to $\sum\limits_{i=1}^n\mu_ip_i+2x\in |\Ocal_E(3)|$. Let $J''=\left\{1\le j\le r:\nu_j''\ne 0\right\}$, consider its complement $J=\{1,2,\dots,r\}\setminus J''$. Set $r''=|J''|\le 3$ and let us analyse the following subcases. 
\vspace{0.5em}

\begin{enumerate}[label=$(\roman*)$]
    \item Take first $\mu_n=1$ and $\nu\ne (5)$ as in (\ref{nu with sum 5 in genus 2}). It is possible to choose $\nu''$ such that it contains at least one entry equal to $1$, let that be $\nu_k''=1$. In particular, $r''\ge 2$ and the points $q_j$ with $j\in J''\setminus\{k\}$ determine $q_k$ uniquely, as the third intersection point of the line through the $q_j$ with $E$. In a similar manner, the points $(p_i)_{i=1}^{n-1}$ determine $p_n$ uniquely via the relation $\sum\limits_{i=1}^n\mu_ip_i+2x\in |\Ocal_E(3)|$, since $\mu_n=1$. Therefore, $\Ycal_2(\mu,\nu)$ is birationally isomorphic to 
$$\left\{\left(E,\ell,(p_i)_{i=1}^{n-1},x,(q_j)_{j\ne k}\right):\begin{matrix}
    E\in |\Ocal_{\PP^2}(3)|\textnormal{, }p_i,x,q_j\in E\textnormal{ for }j\in J''\setminus \{k\}\\
    E\cap \ell\ge \sum\limits_{j\in J''\setminus\{k\}}\nu_j''q_j\textnormal{ and }q_j\in \ell\textnormal{ for }j\in J
\end{matrix}\right\}$$
which, since $n+r''-1\le 9$, is clearly a rational variety of dimension $8+n+r$. 
    \item Take now $\mu_n=1$ with $\nu=(5)$, so $r=r''=1$. For the purposes of this subcase, we fix three points $p_1, q, x\in \PP^2$ using the $\PGL(3)$--action.  If $q$ denotes the (unique) $q_j$, then $q$ is a flex point of $E$ and the embedding of $E$ is given by $|\Ocal_E(3q)|$. Since, as above, the condition $\mu_n=1$ determines $p_n$ uniquely in terms of $(p_i)_{i=1}^{n-1}$, we conclude that $\Ycal_2(\mu,\nu)$ is birationally isomorphic to 
$$\left\{\left(E,(p_i)_{i=1}^{n-1},q,x\right):E\in |\Ocal_{\PP^2}(3)|,\textnormal{ }p_i,x\in E\textnormal{ and }q\in E\textnormal{ is flex}\right\},$$
which in turn is rational for $n \le 7$ by a similar reason to \cite{bellorouskyPHD}. Moreover, we have $\dim \Ycal_2(\mu,\nu)=9+n=8+n+r$. 
    \item The remaining options for $\mu$ are $\mu\in \{(2,2,3),(2,5),(3,4),(7)\}$. We show the variety 
$$\Zcal_2(\mu):=\left\{\left(E,(p_i)_{i=1}^n\right):E\in|\Ocal_{\PP^2}(3)|\textnormal{, }\sum\limits_{i=1}^n\mu_ip_i+2x\in |\Ocal_E(3)|\right\}$$
is rational for these $\mu$, and then, using the same tricks for $\nu$ as in (\ref{nu with sum 5 in genus 2}) above, one concludes easily that $\Ycal_2(\mu,\nu)$ is rational. The idea is just as in the last part of case \textbf{(a)} in the proof of Lemma \ref{rat check Ycal5}: for $\mu=(2,2,3)$, there is a unique point $y\in E$ such that $p_1+p_2+p_3+x\in \Ocal_E(1)(y)$, and then the condition $2p_1+2p_2+3p_3+2x\in |\Ocal_E(3)|$ rewrites as $p_3+2y\in |\Ocal_E(1)|$. Now there are $2$ equations and it is vivid that the data $(E,p_1,p_2,y)$ determines uniquely $p_3,x$. This gives a birational isomorphism 
$$\Zcal_2(\mu)\stackrel{\sim}{\dashrightarrow}\{(E,p_1,p_2,y):E\in|\Ocal_{\PP^2}(3)|\textnormal{, }p_1,p_2,y\in E\},$$
whence $\Zcal_2(\mu)$ is rational. Similar arguments work for the $3$ other options of $\mu$. 
\end{enumerate}

\textbf{\underline{Case $s=6$}.} Let us denote the points with $\nu''_j\ne 0$ by $q^+$ and $q^-$ and define $J:=\{1\le j\le r:q_j\ne q^{\pm}\}$. We split again into subcases.
\vspace{0.5em}

\begin{enumerate}[label=$\roman*)$]
    \item If $\nu$ is not even, it can be directly checked at (\ref{nu with sum 6 in genus 2}) that $\nu''$ can be chosen such that $q^+\ne q^-$. In this situation, $\PP\Omega_2(\mu,\nu)$ is irreducible, and we expect the same from $\Ycal_2(\mu,\nu)$. Indeed, we have
$$\Ycal_2(\mu,\nu)=\left\{\left(E,(p_i)_{i=1}^n,(q_j)_{j\in J},q^{\pm},x\right):\begin{matrix}
    E\in |\Ocal_{\PP^2}(3)|\textnormal{, }p_i,q^+,x\in E, \\
    q_j\in \overline{q^+q^-}\textnormal{ with } q^- \text{ given by}\\
    \sum\limits_{i=1}^n\mu_ip_i+q^{+}+q^-+2x\in |\Ocal_E(4)|\\
\end{matrix}\right\}$$
and this is rational and irreducible precisely when $n+2\le 9$ i.e. $n\le 7$. Furthermore, it is clear that $\dim \Ycal_2(\mu,\nu)=9+n+r$.
    \item Let us continue with the situation $\nu\in \{(2,2,2),(2,4),(6)\}$, when $q^+=q^-=:q$ must hold, but $\mu$ is not even. Again, $\PP\Omega_2(\mu,\nu)$ is irreducible and the same is true for $\Ycal_2(\mu,\nu)$. If $\mu_n=1$, then there is a birational isomorphism
$$\Ycal_2(\mu,\nu)\stackrel{\sim}{\dashrightarrow}\left\{\left(E,(p_i)_{i=1}^{n-1},(q_j)_{j\in J},q,x\right):\begin{matrix}
    E\in |\Ocal_{\PP^2}(3)|\textnormal{, }p_i,q,x\in E\\
    \textnormal{and }\overline{q_jq} \text{ tangent to $E$ at $q$}
\end{matrix}\right\}$$
and the latter variety is rational and irreducible for any partition $\mu$ with $\mu_n=1$, since $n\le |\mu|=8$. This gives $\dim \Ycal_2(\mu,\nu)=9+n+r$ immediately. Finally, the cases $\mu\in \{(2,3,3),(3,5)\}$ are exact copies of the situation that we discussed in genus $5$, regarding the partitions $\mu\in \{(2,2,3,3),(2,3,5)\}$, and thus we omit them. 
    \item Therefore, we are left with the subcases when $\mu$ and $\nu$ are both even. Keep the notation $q^+=q^-=:q$. The condition $\sum\limits_{i=1}^n\mu_ip_i+2q+2x\in |\Ocal_E(4)|$ implies that there exists a line bundle $\tau$ on $E$ such that $\tau^2\cong \Ocal_E$ and 
$$\sum_{i=1}^n\mu_i/2\cdot p_i+q+x\in |\Ocal_E(2)\otimes \tau|.$$
We check that the case $\tau\cong \Ocal_E$ corresponds to the component $\PP\Omega_2(\mu,\nu)^-$. Indeed, if $\tau\cong \Ocal_E$, then there exists a conic $Q$ such that $E\cdot Q=\sum\limits_{i=1}^n\mu_i/2\cdot p_i+q+x$. This gives $\Gamma\cdot (Q+\ell)\ge \sum\limits_{i=1}^n\mu_i/2\cdot p_i+2x+3q+\sum\limits_{j=1}^r\nu_j'q_j$ and, by comparison to the intersection divisor of $\Gamma$ and $E$, it forces $\sum\limits_{i=1}^n\mu_i/2\cdot p_i-\sum\limits_{j=1}^r\nu_j'q_j-q$ to be an effective divisor on $C$. In other words, $\tau\cong \Ocal_E$ does correspond to $\PP\Omega_2(\mu,\nu)^-$. These define the two components of $\Ycal_2(\mu,\nu)$, namely 
$$\Ycal_2(\mu,\nu)^-=\left\{\left(E,(p_i)_{i=1}^n,(q_j)_{j\in J},x\right):\begin{matrix}
    E\textnormal{ cubic, }q\in E\textnormal{, }\overline{q_jq} \text{ tangent to $E$ at $q$},\\
    \sum\limits_{i=1}^n\mu_i/2\cdot p_i+q+x\in |\Ocal_E(2)|
\end{matrix}\right\}$$
and 
$$\Ycal_2(\mu,\nu)^+=\left\{\left(E,\tau,(p_i)_{i=1}^n,(q_j)_{j\in J},x\right):\begin{matrix}
    E\textnormal{ cubic, }q\in E\textnormal{, }\overline{q_jq} \text{ tangent to $E$ at $q$}, \\
    \sum\limits_{i=1}^n\mu_i/2\cdot p_i+q+x\in |\Ocal_E(2)\otimes \tau|\\
    \textnormal{for some }\tau\not\cong \Ocal_E\textnormal{ with }\tau^2\cong \Ocal_E
\end{matrix}\right\}.$$
Since $q$, $x$ appear with coefficient $1$ in the condition on points in both components, one can use the same trick as in the previous subcases, together with the fact that $\Rcal_{1,n}$ is rational for $n\le 6$, to conclude that both $\Ycal_2(\mu,\nu)^-$ and $\Ycal_2(\mu,\nu)^+$ are rational and, furthermore, are of dimension $9+n+r$.
\end{enumerate}

\textbf{\underline{Case $s=7$}.} Finally, this case is concerned with the pole types $$\nu\in \{(2,2,2,1),(4,2,1),(3,2,2),(5,2)\},$$
as already mentioned in (\ref{nu with sum 7 in genus 2}). Let $q_k$ be the unique $q_j$ such that $\nu_j''=1$ and observe that these $4$ partitions have
$$\nu'+2\nu''\in \{(1,1,1,2),(2,1,2),(3,1,1),(4,1)\},$$
respectively. In particular, $\nu'+2\nu''$ contains at least one entry equal to $1$, which we shall call $q_u$, a point distinct from $q_k$; having $u\ne k$ is the precise reason we focus only on these partitions $\nu$. Now let $J=\{1\le j\le r:j\ne k, j\ne u\}$ be the remaining indices. Going through the conditions on points written in the definition of $\Ycal_2(\mu,\nu)$, we see that $(p_i)_{i=1}^n$ and $x$ define $q_k$ uniquely via $\sum\limits_{i=1}^n\mu_ip_i+q_k+2x\in |\Ocal_E(4)|$, and then $q_u$ is uniquely determined by $(p_i)_{i=1}^n$, $x$, $q$ and $(q_j)_{j\in J}$ via the second condition $\sum\limits_{i=1}^n\mu_ip_i+q_u+\sum\limits_{j\ne u}(\nu_j'+2\nu_j'')q_j+2x\in |\Ocal_{E\cup \ell}(4)|$. We deduce that there is a birational isomorphism 
$$\Ycal_2(\mu,\nu)\stackrel{\sim}{\dashrightarrow}\left\{\left(E,\ell,(p_i)_{i=1}^n,(q_j)_{j\ne u},x\right):\begin{matrix}
    E\in |\Ocal_{\PP^2}(3)|\textnormal{, }p_i,x\in E\textnormal{, }\ell\textnormal{ line, }q_j\in\ell\\
    \textnormal{with }q_k\textnormal{ given by }\sum\limits_{i=1}^n\mu_ip_i+q_k+2x\in |\Ocal_{E}(4)|
\end{matrix}\right\}$$
and the latter variety is rational of dimension $9+n+r$, as long as $n+1\le 9$ i.e. $n\le 8$. These settle the claims \ref{rat check for Ycal2 when at most 4 distinct poles} and \ref{dim count gen 2 & at most 4 distinct poles} and also conclude the proof of the Lemma. 
\end{proof}

The remaining part of this subsection is dedicated to \ref{case other mu nu, for sum(nu)=5,6}. Recall that this case is concerned with the pole type $\nu=(1,1,1,1,1)$. We use the following model:

\begin{definition}
Consider the incidence variety 
$$\Xcal_2(\mu,\nu)=\left\{(\Gamma,E,(\ell_k)_{k=1}^2,(p_i)_{i=1}^n,(q_j)_{j=1}^5,x,\qbf):\begin{matrix}
    \Gamma\in |\Ocal_{\PP^2}(4)|\textnormal{ nodal at }x\textnormal{, }\ell_k\textnormal{ lines,}\\
    E\in |\Ocal_{\PP^2}(3)|\textnormal{, }\qbf\in \Gamma_2\textnormal{ with}\\
    \Gamma\cdot E=\sum\limits_{i=1}^n\mu_ip_i+q_2+\qbf+2x,\\
    \Gamma\cdot \ell_1=\sum\limits_{j=2}^5q_j\textnormal{, }\Gamma\cdot \ell_2=q_1+q_2+\qbf
\end{matrix}\right\}.$$
\end{definition}

The first lemma below motivates the construction of $\Xcal_2(\mu,\nu)$ and shows it dominates the space $\PP\Omega_2(\mu,\nu)$. The second lemma is concerned with the proof of the rationality of the same variety $\Xcal_2(\mu,\nu)$ and thus establishes the unirationality of $\PP\Omega_2(\mu,\nu)$ in this range. 

\begin{lemma}\label{dominance in gen 2 for 5 distinct poles}
The map $\Xcal_2(\mu,\nu)\dashrightarrow \PP\Omega_2(\mu,\nu)$ defined by 
$$\left[\Gamma,E,(\ell_k)_{k=1}^2,(p_i)_{i=1}^n,(q_j)_{j=1}^5,x,\qbf\right]\mapsto \left[C,(p_i)_{i=1}^n,(q_j)_{j=1}^5\right],$$
is well-defined and dominant, where $\sigma:C\to \Gamma$ is the normalisation map. Moreover, we obtain a birational isomorphism $\Xcal_2(\mu, \nu)\gitquot \PGL(3) \dashrightarrow \PP\Omega_2(\mu, \nu)$. 
\end{lemma}
\begin{proof}
Take a point $\left[\Gamma,E,(\ell_k)_{k=1}^2,(p_i)_{i=1}^n,(q_j)_{j=1}^5,x,\qbf\right]$ of $\Xcal_2(\mu,\nu)$ and let $\sigma^{-1}(x)=\left\{x^{\pm}\right\}$. We have 
$$\omega_C(x^{+}+x^-)=\sigma^*\omega_\Gamma=\sigma^*\Ocal_{\Gamma}(1)=\sigma^*\Ocal_\Gamma(E-\ell_1-\ell_2)=\Ocal_C\left(x^++x^-+\sum_{i=1}^n\mu_ip_i-\sum_{j=1}^5q_j\right)$$
and thus $\sum\limits_{i=1}^n\mu_ip_i-\sum\limits_{j=1}^5q_j\sim K_C$. Therefore, the map is well-defined. To establish the dominance claim, consider a point $\left[C,(p_i)_{i=1}^n,(q_j)_{j=1}^5\right]$ in $\PP\Omega_2(\mu,\nu)$ such that the points $q_j$ are general on $C$. The line bundle $A:=\Ocal_C\left(\sum\limits_{j=2}^5q_j\right)$ is then a $\gfrak^2_4$ that maps the curve $C$ to a plane quartic $\Gamma\in |\Ocal_{\PP^2}(4)|$ with a node at some point $x\in \PP^2$. 

There is a line $\ell_1$ such that $\Gamma\cdot \ell_1=\sum\limits_{j=2}^5q_j$ and, if $\ell_2$ is the line through $q_1$ and $q_2$, then there is a degree $2$ divisor $\qbf\in \Gamma_2$ such that $\Gamma\cdot \ell_1=q_1+q_2+\qbf$. With these in hand, the relation $\sum\limits_{i=1}^n\mu_ip_i-\sum\limits_{j=1}^5q_j\sim K_C$ implies that 
$$\sum_{i=1}^n\mu_ip_i+q_2+\qbf+x^{+}+x^-\sim K_C\left(\ell_1+\ell_2+x^{+}+x^-\right)\cong \sigma^*\Ocal_\Gamma(3),$$
so there exists a unique cubic $E$ satisfying $\Gamma\cdot E=\sum\limits_{i=1}^n\mu_ip_i+q_2+\qbf+2x$. We conclude that the map is also dominant.
\end{proof}

We are left to prove that the variety $\Xcal_2(\mu,\nu)$ is rational. 

\begin{lemma}
The variety $\Xcal_2(\mu,\nu)$ is rational, for all partitions $\mu$ of $7$. In particular, $\PP\Omega_2(\mu,\nu)$ is unirational, whenever $\nu=(1,1,1,1,1)$ and $\mu$ is a partition of $7$. 
\end{lemma}
\begin{proof}
As in the previous lemmas, let us introduce an auxiliary variety. Define 
$$\Ycal_2(\mu,\nu)=\left\{\left(E,(\ell_k)_{k=1}^2,(p_i)_{i=1}^n,(q_j)_{j=2}^5,\qbf,x\right):\begin{matrix}
    E\in |\Ocal_{\PP^2}(3)|\textnormal{, }\ell_k\textnormal{ lines, }\qbf\in E_2,\\
    E\cdot \ell_2=q_2+\qbf\textnormal{, }\ell_1\cdot \ell_2=q_2\textnormal{, }E\cdot \ell_1\ge q_2,\\
    q_3,q_4,q_5\in \ell_1\textnormal{ and }(\ref{eq: two cond})\textnormal{ holds}
\end{matrix}\right\}$$
and consider the map $\pi:\Xcal_2(\mu,\nu)\dashrightarrow \Ycal_2(\mu,\nu)$ that forgets the curve $\Gamma$ and the point $q_1$. In the definition of $\Ycal_2(\mu, \nu)$ we have made use of the following two conditions:
\begin{align}\label{eq: two cond}
    \sum\limits_{i=1}^n\mu_ip_i+2x\in |\Ocal_E(3)|\textnormal{ and }\sum\limits_{i=1}^n\mu_ip_i+q_3+q_4+q_5+2x\in |\Ocal_{E\cup \ell_1}(3)|.
\end{align}
The next three claims are sufficient to conclude that $\Xcal_2(\mu,\nu)$ is a rational variety:
\begin{enumerate}[label=\textbf{(\roman*)}]
    \item\label{indep lemma gen 2 & 5 distinct poles} the map $\pi$ is a generically injective, 
    \item\label{rat check for Ycal2 when 5 distinct poles} $\Ycal_2(\mu,\nu)$ is an irreducible rational variety for $\mu$ and $\nu$ as in \ref{case other mu nu, for sum(nu)=5,6}, and 
    \item\label{dim count gen 2 & 5 distinct poles} $\dim \Xcal_2(\mu,\nu)=\dim \Ycal_2(\mu,\nu)=14+n$, so $\pi$ is, in fact, a birational isomorphism. 
\end{enumerate}
As before, we use Theorem \ref{big indep lemma} for the statement at $\ref{indep lemma gen 2 & 5 distinct poles}$. Let $S=\textnormal{Bl}_x\left(\PP^2\right)$ be the blow-up of the plane at $x$, write $C\in |4H-2E_x|$ for the strict transform of $\Gamma$ and keep the same notation for $\ell_1$ and $\ell_2$. Then the divisor of conditions is $Z=\sum\limits_{i=1}^n\mu_ip_i+\sum\limits_{j=2}^5q_j+\qbf$, and on $C$ it equals $Z=K_C+2H\vert_C+q_3+q_4+q_5$, so the number of independent conditions imposed on $|4H-2E_x|\cong\PP^{11}$ is $\length(Z)-h^0(C,Z-K_C-2H\vert_C)=11$. This shows that for a generic point $\left[E,(\ell_k)_{k=1}^2,(p_i)_{i=1}^n,(q_j)_{j=2}^5,\qbf,x\right]\in \pi\left(\Xcal_2(\mu,\nu)\right)$, there exists a unique curve $\Gamma$ with a unique singularity as a node at $x$ and satisfying the conditions $\Gamma\cdot E=\sum\limits_{i=1}^n\mu_ip_i+q_2+\qbf+2x$, $\Gamma\cdot \ell_1=\sum\limits_{j=2}^5q_j$ and $\Gamma\cdot \ell_2\ge q_2+\qbf$. For such $\Gamma$, the point $q_1$ is uniquely recovered as the remaining intersection data of $\Gamma$ and $\ell_2$. Therefore, $\pi$ is indeed generically injective. 
\vspace{0.5em}

From Lemma \ref{dominance in gen 2 for 5 distinct poles}, we infer that there is a birational isomorphism $\Xcal_2(\mu,\nu)\gitquot\PGL(3)\dashrightarrow \PP\Omega_2(\mu,\nu)$, and thus $\dim \Xcal_2(\mu,\nu)=14+n$. To ensure $\dim \Ycal_2(\mu,\nu)=14+n$ as well, it suffices to the check the rationality statement \ref{rat check for Ycal2 when 5 distinct poles}, as then the dimensionality statement will follow as a byproduct. So let us focus on the analysis of \ref{rat check for Ycal2 when 5 distinct poles}. We split into subcases. 
\vspace{0.1em}

\begin{enumerate}[label=$\roman*)$]
    \item If $\mu_n=1$, then we remark that the point $p_n$ is recovered uniquely from the curve $E\in |\Ocal_{\PP^2}(3)|$, the points $(p_i)_{i=1}^{n-1}$ and $x\in E$, via the condition $\sum\limits_{i=1}^n\mu_ip_i+2x\in |\Ocal_E(3)|$. Furthermore, the point $q_5$ is also uniquely determined, once the additional points $q_3$ and $q_4$ are given, since $\sum\limits_{i=1}^n\mu_ip_i+q_3+q_4+q_5+2x\in |\Ocal_{E\cup \ell_1}(3)|$ must hold. Therefore, $\Ycal_2(\mu,\nu)$ is birationally isomorphic to 
$$\left\{\left(E,(\ell_k)_{k=1}^2,(p_i)_{i=1}^{n-1},(q_j)_{j=2}^4,x\right):\begin{matrix}
    E\in |\Ocal_{\PP^2}(3)|\textnormal{, }\ell_k\textnormal{ lines, }p_i,x\in E,\\
    q_2\in E\cap \ell_1\cap \ell_2\textnormal{, }q_3,q_4\in \ell_1
\end{matrix}\right\}.$$
Since $n\le |\mu|=7$, it follows that the latter variety is rational and one can easily see now that $\dim \Ycal_2(\mu,\nu)=14+n$, as required.
    \item The remaining partitions to check are $\mu\in\{(2,2,3),(2,5),(3,4),(7)\}$. In any of these cases, a combination of the argument in the subcase $\mu_n=1$ just above and the argument regarding the rationality of $\Zcal_2(\mu)$ presented in the case $s=5$ of part \ref{case of at most 4 distinct poles} will do. 
\end{enumerate}

\end{proof}

\bibliography{bibliography}
\bibliographystyle{amsalpha}
\end{document}